\documentclass[matprg,final]{svjour}

\usepackage{amsmath}
\usepackage{amsfonts}
\usepackage{amssymb}
\usepackage{graphics}
\usepackage{epsfig}
\usepackage{url}

\newtheorem{lem}[theorem]{Lemma}

\newtheorem{cor}[theorem]{Corollary}
\newtheorem{prop}[theorem]{Proposition}

\newtheorem{conj}[theorem]{Conjecture}

\newtheorem{exmp}[theorem]{Example}
\newtheorem{exm}[theorem]{Example}

\newcommand{\CC}{\mathbb C}
\newcommand{\RR}{\mathbb R}
\newcommand{\QQ}{\mathbb Q}
\newcommand{\Q}{\mathbb{Q}}
\newcommand{\ZZ}{\mathbb Z}

\newcommand{\PP}{\mathbb P}

\newcommand{\Oo}{\mathcal O}
\newcommand{\U}{\mathcal U}
\newcommand{\V}{\mathcal V}
\newcommand{\W}{\mathcal W}

\newcommand{\reff}[1]{(\ref{#1})}
\newcommand{\mc}[1]{\mathcal{#1}}

\newcommand{\bbm}{\begin{bmatrix}}
\newcommand{\ebm}{\end{bmatrix}}

\begin{document}

\title{The Algebraic Degree of Semidefinite Programming}

\author{Jiawang Nie \and Kristian Ranestad \and Bernd Sturmfels}

\institute{
Jiawang Nie \at
Department of Mathematics, UC San Diego,
La Jolla, CA 92093, USA.
\email{njw@math.ucsd.edu} \and
Kristian Ranestad \at
Department of Mathematics, University of Oslo,
PB 1053 Blindern, 0316 Oslo, Norway.
\email{ranestad@math.uio.no} \and
Bernd Sturmfels \at
Department of Mathematics, UC Berkeley,
   Berkeley, CA 94720, USA.
\email{bernd@math.berkeley.edu}
}

\maketitle

\begin{abstract}
Given a generic semidefinite program, specified by matrices
with rational entries, each coordinate of its optimal solution is
an algebraic number.  We study the degree of
the minimal polynomials of these algebraic numbers. Geometrically,
this degree counts the critical points attained by a linear functional
on a fixed rank locus in a linear space of
symmetric matrices. We determine this degree
using methods from complex algebraic geometry, such
as projective duality, determinantal varieties, and their Chern classes.
\end{abstract}


\maketitle

\section{Introduction}

A fundamental question about any optimization problem is
how the output depends on the input. The set of optimal solutions
and the optimal value of the problem are functions of the parameters,
and it is important to understand the nature of these functions.
For instance, for a linear programming problem,
\begin{equation}
\label{LP} {\rm maximize} \,\,\, c \cdot x \,\,\,
{\rm subject} \, {\rm to} \,\,\,A \cdot x = b \,\,\,{\rm and} \,\,\,x \geq 0 ,
\end{equation}
the optimal value is convex and piecewise linear in
the cost vector $c$ and the right hand side $b$,
and it is a piecewise rational function of the entries
of the matrix $A$. The area of mathematics which studies
these functions is geometric combinatorics, specifically
the theory of matroids for the dependence on $A$,
and the theory of polyhedral subdivisions \cite{DRS}
for the dependence on $b$ and $c$.

For a second example, consider the following basic question in game theory:
\begin{equation}
\label{Nash}
\hbox{Given a game, compute its Nash equilibria.}
\end{equation}
If there are only two players and
one is interested in fully mixed Nash equilibria
then this is a linear problem, and in fact closely
related to linear programming.
On the other hand, if the number of players is more
than two  then the problem (\ref{Nash}) is
universal in the sense of real algebraic geometry:
Datta \cite{Dat} showed that every real algebraic
variety is isomorphic to the  set of Nash equilibria
of some three-person game. A corollary of her
construction is that, if the Nash equilibria are
discrete, then their coordinates can be
arbitrary algebraic functions of the given
input data (the payoff values which specify the game).

Our third example concerns maximum likelihood
estimation in statistical models for discrete data.
Here the optimization problem is as follows:
\begin{equation}
\label{MLE1}
{\rm Maximize} \,\,\,
 p_1(\theta)^{u_1} p_2(\theta)^{u_2} \cdots p_n(\theta)^{u_n}
 \,\,\, {\rm subject}\, \,\, {\rm to} \,\,\, \theta \in \Theta,
 \end{equation}
 where $\Theta$ is an open subset of $\RR^m$,
 the $p_i(\theta)$ are polynomial functions that sum to one,
 and the $u_i$ are positive integers (these are the data).
The optimal solution $\hat \theta$, which is the maximum likelihood
estimator,
depends on the data:
\begin{equation}
\label{MLE2}
 (u_1,\ldots,u_n) \,\mapsto \, \hat{\theta}( u_1,\ldots,u_n) .
 \end{equation}
 This is an algebraic function, and recent work in algebraic
 statistics \cite{CHKS,HKS} has led to a formula for the
 degree of this algebraic function,
 under certain hypothesis on the polynomials $p_i(\theta)$
 which specify the statistical model.

The aim of the present paper is to conduct a similar
algebraic analysis for the optimal value function in
{\em semidefinite programming}.
This function shares some key features with each of the
previous examples. To begin with, it
is a convex function which is piecewise algebraic.
However, unlike for (\ref{LP}), the pieces are non-linear, so
there is a notion of algebraic degree as for
Nash equilibria (\ref{Nash}) and
maximum likelihood estimation (\ref{MLE1}).
However, semidefinite programming
does not exhibit universality as in \cite{Dat}
because the structure of  real symmetric matrices
imposes some serious constraints on the
underlying algebraic geometry. It is these
constraints we wish to explore and characterize.

We consider the semidefinite programming (SDP) problem in the form
\begin{equation}
\label{SDP} \,\,{\rm maximize} \,\,\,{\rm trace}(B \cdot Y ) \,\,\,
\hbox{subject to}\,\,\, Y \in \mathcal{U}  \,\, \,{\rm and} \,\,\, Y \, \succeq \, 0.
\end{equation}
where $B$ is a real symmetric $n \times n$-matrix,
$\mathcal{U}$ is a $m$-dimensional affine subspace
in the $\binom{n+1}{2}$-dimensional space of
real $n \times n$-symmetric matrices,
and $\,Y \,\succeq \, 0 \,$ means that
$Y$ is positive semidefinite
(all $n$ eigenvalues of $Y$ are non-negative).
The problem (\ref{SDP}) is {\em feasible} if and only if
the subspace $\mathcal{U}$ intersects the
cone of positive semidefinite matrices.
In the range where it exists and is unique,
the optimal solution $\hat Y$ of this problem is a piecewise
algebraic function of the matrix $B$ and
the subspace $\mathcal{U}$. Our aim is
to understand the geometry of this function.

Note if  $\mathcal{U}$ consists of diagonal matrices
only, then we are in the linear programming case (\ref{LP}),
and the algebraic degree of the pieces of $\hat Y$ is just one.
What we are interested in is the increase in algebraic
complexity which arises from the passage from diagonal matrices
to non-diagonal symmetric matrices.

\begin{figure}
\begin{center}
\includegraphics[totalheight=7.5cm]{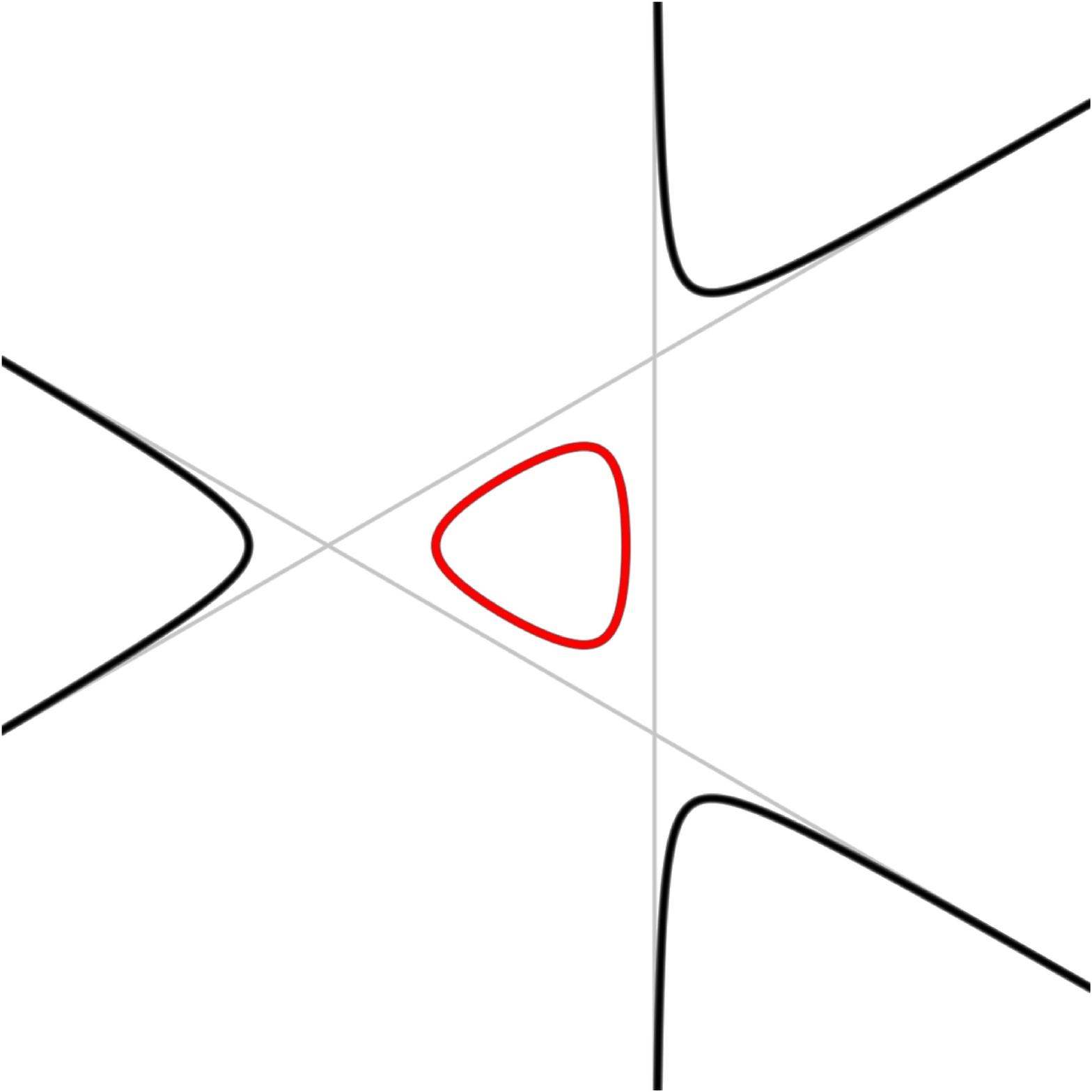}
\end{center}
\caption{The convex component in the center of this
elliptic Vinnikov curve
 is the feasible region for SDP with $m=2,n=3$.}
     \label{fig:vinnikov}  \end{figure}

\begin{figure}
\begin{center}
\includegraphics[totalheight=7.5cm]{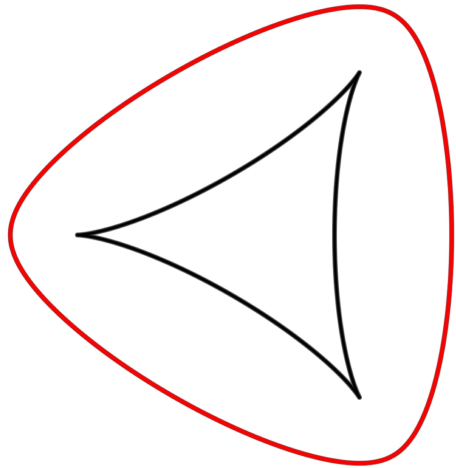}
\end{center}
\vskip -.6cm
\caption{The dual to the elliptic Vinnikov curve in Figure \ref{fig:vinnikov}
is a plane curve of degree six
with three real singular points.}
     \label{fig:vinnidual} \end{figure}

\begin{exm}[Elliptic Vinnikov curves] \label{ex:vinnikov}
Let $n=3$ and $m=2$, so $Y$ runs over a
two-dimensional affine space
of symmetric $3 \times 3$-matrices.
Then $Y \succeq 0$ specifies
a closed convex semi-algebraic region in this
plane, whose boundary is the central
connected component of a cubic curve
as depicted in Figure~\ref{fig:vinnikov}.
This curve is a {\em Vinnikov curve},
which means that it satisfies
the following constraint  in real
algebraic geometry: any line which meets the interior
of the convex region intersects this cubic curve
in three real points.  See \cite{HV,LPR,Spe,Vin} for details.
However, there are no constraints on Vinnikov curves in the setting
 of complex algebraic geometry.
The Vinnikov constraints involve inequalities and
no equations. This explains why the curve in Figure \ref{fig:vinnikov}
is smooth.

   Our problem (\ref{SDP}) is to maximize a linear
function over the convex component. Algebraically,
the restriction of a linear function to the cubic curve
has six critical points, two of which are complex.
They correspond to the intersection points of the
dual curve with a line in the dual projective plane.
The degree six curve dual to the elliptic Vinnikov curve is
depicted in Figure \ref{fig:vinnidual}.

Our analysis shows
that the algebraic degree of SDP equals six when  $m=2$ and $n=3$.
If the matrix $B$ and the plane $\mathcal{U}$
are defined over $\QQ$ then
the coordinates of the optimal solution $\hat Y$ are algebraic
numbers of degree six.
By Galois theory, the solution $\hat Y$ cannot in general
be expressed in terms of radicals. 
For any specific numerical instance we can use the command
``galois" in {\tt maple} to compute the Galois group, which is
then typically found to be the symmetric group $S_6$.
\qed
\end{exm}

\begin{exm}[The Cayley Cubic] \label{ex:cayley}
 Now, suppose that $m=n=3$. Then ${\rm det}(Y) = 0 $ is a cubic surface,
but this surface is constrained
in the context of complex algebraic geometry because it
cannot be smooth. The cubic surface ${\rm det}(Y) = 0$
has four isolated nodal singularities, namely,
the points where $X$ has rank one.
This cubic surface is known to geometers as
the {\em Cayley cubic}. In the optimization literature, it occurs
under the names {\em elliptope} or {\em symmetroid}.
Optimization experts are familiar with (the convex component of)
the Cayley cubic surface from (the upper left hand picture in) Christoph Helmberg's
SDP web page \
{\tt http://www-user.tu-chemnitz.de/$\sim$helmberg/semidef.html}.

\begin{figure}
\begin{center}
\includegraphics[totalheight=8cm]{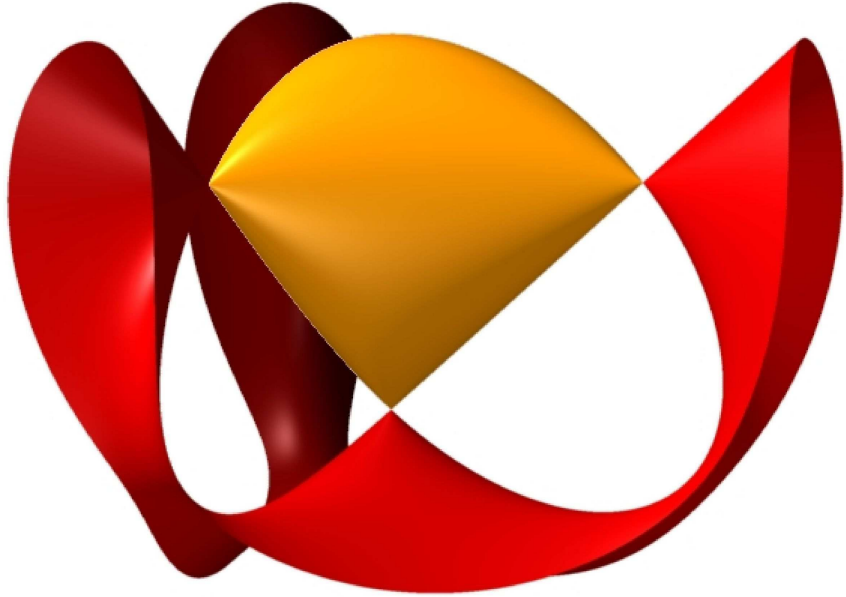}
\end{center}
\caption{The convex component in the center of this
Cayley cubic surface
 is the feasible region for SDP with $m=n=3$.}
     \label{fig:cayley}  \end{figure}

The surface dual to the Cayley cubic is a surface of degree four, which is
known as the {\em Steiner surface}. There
are now two possibilities for the optimal solution $\hat Y$ of
(\ref{SDP}). Either $\hat Y$ has rank one, in which case it is one of the
four singular points of the cubic surface in Figure \ref{fig:cayley},
or $\hat Y$ has rank two
and is gotten by intersecting the Steiner surface
by a line specified by $B$. In either of these two cases,
the optimal solution $\hat Y$ is an algebraic function of degree
four in the data specifying $B$ and $\mathcal{U}$.
In particular, using Girolamo Cardano's {\it Ars Magna},
 we can express the coordinates of $\hat Y$
in terms of radicals in $(B,\U)$. \qed
\end{exm}

 The objective of this paper is to study the geometric figures
shown in Figures 1, 2 and 3 for
arbitrary values of $n$ and $m$. The targeted audience
includes both algebraic geometers and scholars in optimization.
Our presentation is organized as follows. In Section 2 we
review SDP duality, we give an elementary introduction to
the notion of algebraic degree,
and we explain what it means for the data
$\mathcal{U}$ and $B$ to be generic.
In Section 3 we derive Pataki's inequalities which
characterize the possible ranks of the
optimal matrix $\hat Y$, and, in Theorem \ref{zdim},
we present a precise characterization of
the algebraic degree. The resulting geometric
formulation of semidefinite programming is our
point of departure in Section 4. Theorem
\ref{bideg} expresses the algebraic degree of SDP
as a certain bidegree. This is a notion of degree for
subvarieties of products of projective spaces,
and is an instance of the general definition
of multidegree  in Section 8.5 of the text book \cite{MS}.
Theorem \ref{degformula} gives
explicit formulas for the algebraic degree,
organized according to the rows of Table \ref{degtab}.
In Section 5 we present results involving
projective duality and determinantal varieties.
in Section 6 this is combined with results of Pragacz \cite{Prag}
to prove Theorem~\ref{degformula}, and to derive the general
 formula stated in Theorem~\ref{smooth} and Conjecture~\ref{conj}.

\smallskip

Two decades ago, the concept of algebraic degree of an
optimization problem had been explored in the
computational geometry literature, notably in the work
of Bajaj \cite{Ba}. However, that line of research had only
few applications, possibly because of the dearth of precise results
for geometric problems of interest. Our paper fills this gap,
at least for problems with semidefinite representation, and it
can be read as an invitation to experts in complexity theory
to take a fresh look at Bajaj's conclusion  that
{\em ``... the domain of relations between the
algebraic degree ... and the complexity of obtaining
the solution point of optimization problems
is an exciting area to explore.''} \cite[page 190]{Ba}.

\smallskip

The algebraic degree of semidefinite programming
addresses the computational complexity at a fundamental level.
To solve the semidefinite programming exactly
essentially reduces to solve
a class of univariate polynomial equations
whose degrees are the algebraic degree.
As we will see later in this paper,
the algebraic degree is usually very big,
even for some small problems.
An explicit general formula for the algebraic degree
was given by von Bothmer and Ranestad
in the paper \cite{vBR} which was written subsequently
to this article.

\section{Semidefinite programming: duality and symbolic solutions}
\setcounter{equation}{0}

In this section we review the duality theory of
semidefinite programming, and we give an elementary introduction to the
notion of algebraic degree.

Let $\RR$ be the field of real numbers and $\Q$ the subfield
of rational numbers. We write ${\mathcal S}^n$ for the
$\binom{n+1}{2}$-dimensional vector space of
symmetric $n \times n$-matrices over $\RR$,
and $\Q{\mathcal S}^n$ when we only allow entries in $\Q$.
A matrix $X \in \mathcal{S}^n$ is {\em positive definite},
denoted $X \succ 0$, if  $u^T X u >  0$ for all $u\in \RR^n \backslash \{0\}$,
and it is {\em positive semidefinite}, denoted $X \succeq 0$,
if $u^T X u \geq 0$ for all $u \in \RR^n$.
We consider the {\it semidefinite programming} (SDP) problem
\begin{align}
\min_{X\in {\mathcal S}^n} & \quad  C \bullet X \label{sdp1}  \\
\hbox{s.t.} & \quad A_i \bullet X =b_i \quad \hbox{for} \,\,\, i\, =1,\ldots,m \label{sdp2}\\
\hbox{and} & \quad X \succeq 0\label{sdp3}
\end{align}
where $b\in \Q^m, C,\,A_1,\ldots,A_m \in \Q\mathcal{S}^n$.
The inner product $C\bullet X$
is defined as
$$ C \bullet X \,\,\, := \,\,\, {\rm trace}(C \cdot X)
\,\,= \,\,\sum C_{ij} X_{ij} \qquad
\hbox{for  $C,X\in {\mathcal S}^n$}. $$
This is a linear function in $X$ for fixed $C$,
and that is our {\em objective function}.
The primal SDP problem \reff{sdp2}-\reff{sdp3} is called
  {\em strictly feasible}, or the feasible region
{\em has an interior point}, if there exists some $X\succ 0$
such that \reff{sdp2} is met.

Throughout this paper, the words ``generic'' and
``genericity'' appear frequently. These notions have a precise meaning in
algebraic geometry: the data $C,b,A_1,\ldots,A_m$ are  {\em generic} if
their coordinates satisfy no non-zero polynomial
equation with coefficients in $\QQ$. Any statement that is
proved under such a genericity hypothesis will be valid for
all data that lie in a dense, open subset of the space of data,
and hence it will hold except on a set of Lebesgue measure zero. For a simple
illustration consider the quadratic equation $\alpha t^2+\beta t+\gamma=0$
where $t$ is the variable and
 $\alpha,\beta,\gamma$ are certain parameters. This equation has two distinct
roots if and only if the discriminant $\alpha(\beta^2-4\alpha\gamma)$ is non-zero.
The equation $\alpha(\beta^2-4\alpha\gamma )= 0$ defines a surface,
which has measure zero in $3$-space.
The general point $(\alpha,\beta,\gamma)$ does not lie on this surface.
So we can say that $\alpha t^2+\beta t+\gamma=0$ has two distinct
roots when $\alpha,\beta,\gamma$ are generic.

The convex optimization problem dual to \reff{sdp1}-\reff{sdp2} is as follows:
\begin{align}
\max_{y\in \RR^m} & \quad  b^Ty  \label{lmi1} \\
{\rm s.t.} & \quad  A(y) \,\, := \,\, C-\sum_{i=1}^m y_iA_i \quad \succeq \quad 0. \label{lmi2}
\end{align}
Here the decision variables $y_1,\ldots,y_m$ are  real unknown.
The condition \reff{lmi2} is also called a {\it linear matrix inequality} (LMI).
We say that $\reff{lmi2}$ is {\em strictly feasible} or
has an {\em interior point} if there exists some $y\in \RR^m$
such that $A(y)\succ 0$.

Our formulation of semidefinite programming in (\ref{SDP})
is equivalent to \reff{lmi1}--\reff{lmi2} under the following identifications.
Take $\mathcal{U}$ to be the affine space consisting of all matrices
$\,C-\sum_{i=1}^m y_iA_i \,$ where $ y \in \RR^m$,
write $Y = A(y)$ for an unknown matrix in this space, and fix a matrix $B$
such that $B \bullet A_i = -b_i$ for $i = 1,\ldots,m$.
Such a choice is possible provided the
matrices $A_i$ are chosen to be linearly independent, and it implies
 $\,B \bullet Y - B \bullet C \,=\, b^T y$.

We refer to \cite{VB,WSV} for the theory, algorithms and applications of SDP.
The known results concerning SDP duality can be summarized as follows.
Suppose that $X \in {\mathcal S}^n$ is feasible for \reff{sdp2}-\reff{sdp3}
and $y \in \RR^m$ is feasible for \reff{lmi2}. Then
$\, C\bullet X - b^T y \, = \, A(y)\bullet X \,\geq \,0$,
because the  inner product of any two semidefinite matrices is non-negative.
Hence the following {\it weak duality} always holds:
\[
\sup_{A(y)\succeq 0} b^Ty \quad \leq \,\,\,\, \inf_{
\substack{
 X\succeq 0 \\
\forall i \,:\,
A_i \bullet X =b_i}}
C\bullet X .
\]
When equality holds in this inequality, then we say
that {\it strong duality} holds.

The cone of positive semidefinite matrices is a self-dual cone.
This has the following important consequence for any
positive semidefinite matrices $A(y)$ and $X$ which
represent feasible solutions for \reff{sdp2}-\reff{sdp3}  and
 \reff{lmi2}.  The inner  product  $\,A(y)\bullet\,X\,$ is
 zero  if and only if the
 matrix product $\,A(y) \cdot X\,$ is the
 zero matrix.
The optimality conditions are summarized
in the following theorem:

\begin{theorem}[Section 3 in \cite{VB} or
Chapter 4 in \cite{WSV}] \label{optcnd}
\label{wsv}
Suppose that both the primal problem
 \reff{sdp1}-\reff{sdp3} and the dual problem
\reff{lmi1}-\reff{lmi2} are strictly feasible.
Then strong duality holds,
there exists a pair of optimal solutions,
and the
following optimality conditions characterize
the pair of optimal solutions:
\begin{align}
A_i \bullet \hat X & =b_i \quad {\rm for} \,\,\,\, i =1,2,\ldots,m, \label{optcond1} \\
A(\hat y) \cdot \hat X & = 0,  \label{optcond2} \\
A(\hat y) \succeq 0 \,\,\, {\rm and} \,\,\, \hat X &\succeq 0. \label{optcond3}
\end{align}
\end{theorem}
The matrix equation~\reff{optcond2} is the {\it
complementarity condition}. It implies that the sum of the
ranks of the matrices $A(\hat y)$ and $\hat X$ is at most $n$.
We say that {\it strict complementarity} holds if
the sum of the ranks of $A(\hat y)$ and
$\hat X$ equals $n$.

Suppose now that the given data $C,A_1,\ldots,A_m$ and $b$
are generic over the rational numbers $\QQ$. In practice,
we may choose the entries of these matrices to be
random integers, and we compute
the optimal solutions $\hat y$ and $\hat X$ from these data,
using a numerical interior point method.
Our objective is to learn by means of algebra how
$\hat y$ and $\hat X$ depend on the input
$C,A_1,\ldots,A_m, b$. Our approach
rests on the optimality conditions
in Theorem \ref{optcnd}. These take the form of
a system of polynomial equations. If we could solve these
equations using symbolic computation, then this would furnish an exact
representation of the optimal solution $(\hat X, \hat y)$.
We illustrate this approach for a small example.

\begin{exmp} \label{numericsymbolic}
Consider the following semidefinite programming problem:
\begin{align*}
{\rm Maximize} & \quad y_1+y_2+y_3 \\
{\rm subject} \,\,{\rm to} & \,\,\, A(y) \, := \,
\bbm
y_3+1 & y_1+y_2 & y_2 & y_2+y_3 \\
y_1+y_2 & -y_1+1 & y_2-y_3 &y_2 \\
y_2 &  y_2-y_3 & y_2+1 & y_1+y_3 \\
y_2+y_3 & y_2 &   y_1+y_3 & -y_3+1
\ebm \,\,\succeq \,\,0.
\end{align*}
This is an instance of the LMI formulation \reff{lmi1}-\reff{lmi2} with $m=3$ and $n=4$.
Using the numerical software {\tt SeDuMi} \cite{sedumi}, we
easily find the  optimal solution:
$$ (\hat y_1, \hat y_2, \hat y_3) \quad = \quad
\bigl(\,
      0.3377\dots,\,
      0.5724\dots,\,
      0.3254\dots \bigr). $$
What we are interested in is to understand the nature of these three numbers.

To examine this using symbolic computation,
we also consider the primal formulation \reff{sdp1}--\reff{sdp3}.
Here the decision variable is a symmetric
$4 \times 4$-matrix $X = (x_{ij})$ whose
ten entries $x_{ij}$ satisfy the three linear constraints \reff{optcond1}:
$$ \begin{matrix}
 A_1 \bullet X \,\,\,   = & \,\,\,   -2 x_{12}+x_{22}-2 x_{34}
& = \quad 1 \\
 A_2 \bullet X   \,\,\, = & \,\,\, -2 x_{12} - 2 x_{13} - 2 x_{14}-2 x_{23}-2 x_{24}-x_{33}
& = \quad 1 \\
 A_3 \bullet X  \,\,\, = & \,\,\, x_{11} - 2 x_{14} + 2 x_{23} - 2 x_{34} + x_{44}
& = \quad 1
\end{matrix} $$
In addition, there are sixteen quadratic constraints coming from the
complementarity condition \reff{optcond2}, namely we
equate each entry of the matrix product
$\,A(y) \cdot X \,$ with zero. Thus \reff{optcond1}--\reff{optcond2}
translates into a system
of $19$ linear and quadratic polynomial equations in the $13$ unknowns
$\,y_1, y_2, y_3 ,x_{11}, x_{12}, \ldots, x_{44}$.

Using symbolic computation methods (namely, Gr\"obner bases
in {\tt Macaulay 2} \cite{M2}), we find that
these equations have finitely many complex solutions.
The number of solutions is $26$. Indeed, by
eliminating variables, we discover that each coordinate $y_i$
or $x_{jk}$ satisfies a univariate equation of degree $26$.
Interestingly, these univariate polynomials are not
irreducible but they factor. For instance, the optimal
first coordinate ${\hat y}_1$ satisfies the univariate
polynomial $f(y_1)$ which factors into
a polynomial $g(y_1)$ of degree $16$ and a polynomial
of degree $h(y_1)$ of degree $10$. Namely, we have
$\, f(y_1) \,\, = \,\, g(y_1) \cdot h(y_1) $, where
\begin{align*}
g(y_1) \,=
&\quad 403538653715069011\, y_1^{16}-2480774864948860304\, y_1^{15} \\ &
+6231483282173647552\, y_1^{14}-5986611777955575152\, y_1^{13} \\ &
 + \qquad \cdots \qquad \cdots \qquad \cdots  \qquad \cdots \qquad + \\ &
+59396088648011456\, y_1^{2}-4451473629111296\, y_1+149571632340416
\end{align*}
and $\, h(y_1) \,\,= \,\,
2018\, y_1^{10}-12156\, y_1^{9}+17811\, y_1^{8}+ \, \cdots \,+1669\,y_1-163 $.

Both of these polynomials are irreducible in $\QQ[y_1]$.
By plugging in, we see that the optimal first coordinate
$\,\hat y_1 = 0.3377\dots\,$ satisfies $g(\hat y_1) = 0$.
Hence $\hat y_1$ is an algebraic number  of degree $16$ over $\QQ$.
Indeed, each of the other twelve optimal coordinates
$\,\hat y_2, \hat y_3, \hat x_{11}, \hat x_{12}, \ldots, \hat x_{44}\,$
also has degree $16$ over $\QQ$.
We conclude that the algebraic degree of this particular SDP problem is $16$.
Note that the optimal matrix $A(\hat y)$ has rank $3$ and the matrix $\hat X$ has rank $1$.

We are now in a position to vary the input data and
perform a parametric analysis. For instance, if the
objective function is changed to $y_1-y_2$
then the algebraic degree is $10$ and the
ranks of the optimal matrices are both $2$. \qed
\end{exmp}

The above example is not special. The entries for $m=3,n=4$ in
Table~\ref{degtab} below inform us that, for generic data $(C,b,A_i)$,
the algebraic degree is
 $16$ when the optimal matrix $\hat Y = A(\hat y)$ has rank three,
the algebraic degree is $10$ when the optimal matrix $\hat Y  = A(\hat y)$
has rank two, and rank one or four optimal solutions do not exist.
These former two cases can be understood geometrically by drawing a picture
as in Figure 3 and Example \ref{ex:cayley}. The surface ${\rm det}(Y) = 0$
has degree four and it has $10$ isolated singular points (these are
 the matrices $Y$ of rank two). The surface dual to the
quartic ${\rm det}(Y) = 0$ has degree $16$. Our optimal solution
in Example \ref{numericsymbolic} is one of the $16$ intersection points
of this dual surface with the line specified by the linear objective function
$\,B \bullet Y = b \cdot y^T$.
The concept of duality in algebraic geometry will be reviewed in Section 5.

For larger semidefinite programming problems it is  impossible
(and hardly desirable)
 to solve the polynomial equations \reff{optcond1}--\reff{optcond2}
symbolically.  However, we know that the coordinates
$\hat y_i$ and $\hat x_{jk}$ of the optimal solution
are the roots of some univariate polynomials
which are irreducible over $\Q$. If the data
are generic, then the degree of these univariate polynomials depends only on the rank of the optimal solution.
This is what we call the {\it algebraic degree}
of the semidefinite programming problem~\reff{sdp1}-\reff{sdp3} and its dual
\reff{lmi1}-\reff{lmi2}. The objective of this paper is to find a formula for this algebraic degree.

\section{From Pataki's inequalities to algebraic geometry}
\setcounter{equation}{0}

Example \ref{numericsymbolic} raises the question which ranks
are to be expected for the optimal matrices.
This question is answered by the following result
from the SDP literature \cite{AHO,Pat}.
 We refer to the semidefinite program
specified in \reff{sdp1}--\reff{sdp3}, \reff{lmi1}--\reff{lmi2}
and \reff{optcond1}--\reff{optcond3}. Furthermore, we always
assume that the problem instance $(C,b,A_1,\ldots,A_m)$
is generic, in the sense discussed above.

\begin{prop}{{\rm (Pataki's Inequalities \cite[Corollary 3.3.4]{Pat}, see
also \cite{AHO})}} \label{pataki}   \hfill \break
Let $r$ and $n-r$ be the ranks
of the optimal matrices $\hat Y  = A(\hat y)$ and $\hat X$. Then
\begin{align} \label{patiq}
 \binom{n-r+1}{2}\, \leq\, m \qquad
 \hbox{and} \qquad \binom{r+1}{2} \,\leq\, \binom{n+1}{2}-m.
\end{align}
\end{prop}

A proof will be furnished later in this section.
First we illustrate Proposition \ref{pataki}
by describing a numerical experiment which is easily performed
for a range of values  $(m, n)$.
We generated $m$-tuples of $n \times n$-matrices $(A_1,A_2, \ldots,A_m)$,
where each entry was independently
drawn from the standard Gaussian distribution.
Then, for a random positive semidefinite matrix $X_0$,
let $b_i= A_i \bullet X_0$,
which makes the feasible set \reff{sdp2}-\reff{sdp3} is nonempty.
For each such choice, we generated $10000$ random
symmetric matrices $C$. Using {\tt SeDuMi} \cite{sedumi} we
then solved the program \reff{sdp1}-\reff{sdp3} and
we determined the numerical rank of its optimal solution $\hat Y$.
This was done by computing the Singular Value Decomposition of $\hat Y$
in {\tt Matlab}. The result is
the rank distribution in Table~\ref{rkdistrib}.

\begin{table} \small
\begin{center}
\begin{tabular}{|c|cr|cr|cr|cr|} \hline
$n$ & \multicolumn{2}{c|}{ 3} & \multicolumn{2}{c|}{4} &
\multicolumn{2}{c|}{5} & \multicolumn{2}{c|}{6} \\ \hline
$m$ & rank & percent& rank & percent
& rank & percent  & rank & percent \\  \hline

3&  2& 24.00\% &  3& 35.34\% &   4& 79.18\% &  5& 82.78\% \\
  &  1& 76.00\% &  2& 64.66\% &   3& 20.82\% &  4& 17.22\% \\  \hline

4&  &  &             3& 23.22\%  &   4&16.96\%   &  5& 37.42\% \\
  &  1& 100.00\%  &  2& 76.78\%  &   3&83.04\%   &  4& 62.58\% \\
\hline

5&  & &  &  &                           4&5.90\%  &   5& 38.42\%   \\
  &  1 &  100.00\% & 2  & 100.00 \% &   3&94.10\%   &  4& 61.58\%    \\
\hline

6 &   &&   &    &     &                  & 5&  1.32\%\\
    &      &&  2&67.24\%  & 3&93.50\%    & 4& 93.36\%  \\
   &       &&  1&32.76\%  &  2&6.50\%    & 3&  5.32\% \\  \hline

7&  &&   2&52.94\%  &    3& 82.64\% &  4& 78.82\%  \\
  &  &&  1&47.06\%  &    2& 17.36\% &  3& 21.18\%  \\  \hline

8&  &&  &  &                3&34.64\%  &  4& 45.62\% \\
  &  &&  1 & 100.00\%  &    2&65.36\%  &  3& 54.38\%  \\  \hline

9& &&  &  &                 3&7.60\%   &  4& 23.50\%   \\
  &  &&    1& 100.00\% &    2&92.40\%  &  3& 76.50\%   \\  \hline
\end{tabular}
\smallskip
\end{center}
\caption{Distribution of the rank of the optimal matrix $\hat Y$.}
 \label{rkdistrib}
\end{table}

Table~\ref{rkdistrib} verifies that, with probability one, the rank $r$ of
the optimal matrix $\hat Y$ lies in the interval specified by Pataki's Inequalities.
The case $m=2$, which concerns Vinnikov curves as in
Example \ref{ex:vinnikov}, is not listed in Table~\ref{rkdistrib}
because here the optimal rank is always $n-1$.
The first interesting case is  $m=n=3$, which concerns
Cayley's cubic surface as in Example \ref{ex:cayley}.
When optimizing a linear function over the convex
surface in Figure~\ref{fig:cayley}, it is three times less likely
for a smooth point to be optimal than one of the four vertices.
For $m=3,n=4$, as in Example~\ref{numericsymbolic},
the  odds are slightly more balanced.
In only $35.34 \%$ of the instances the algebraic degree of
the optimal solution was $16$, while in
$64.66 \%$ of the instances the algebraic degree was
found to be $10$.

This experiment highlights again the genericity hypothesis made
throughout this paper. We shall always assume that
the $m$-tuple $(A_1,\ldots,A_m)$, the cost matrix $C$ and vector $b$ are generic.
All results in this paper are only valid under this
genericity hypothesis. Naturally, special
phenomena will happen for special problem instances.
For instance, the rank of the optimal matrix can be outside
the Pataki interval. While
such special instances are important for applications of SDP,
we shall not address them in this present study.

In what follows we introduce our algebraic setting.
We fix the assumptions
\begin{equation}
\label{Bassumption}
b_1=1 \,\, \,\hbox{and} \,\,\, b_2=b_3=\cdots=b_m=0.
\end{equation}
These assumptions appear to violate our genericity hypothesis.
However, this is not the case since any
generic instance can be transformed,
by a linear change of coordinates, into an instance
of \reff{optcond1}--\reff{optcond3} that
satisfies (\ref{Bassumption}).

Our approach is based on two linear spaces of symmetric matrices,
\begin{align*}
\U &\quad = \quad \langle C,A_1,A_2, \ldots ,A_m\rangle  \,\, \subset \,\,\mathcal{S}^n ,\\
\V &\quad = \qquad  \langle A_2, \ldots , A_m \rangle \qquad \subset \,\, \U.
\end{align*}
Thus, $\U$ is a generic linear subspace of dimension $m+1$ in  $\mathcal{S}^n$,
and $\V$ is a generic subspace of codimension $2$ in $\U$.
This specifies a dual pair of flags:
\begin{equation}
\label{flags}
\V \subset \U \subset \mc{S}^n \quad \hbox{and} \quad
\U^\perp \subset \V^\perp \subset (\mc{S}^n)^*.
\end{equation}
In the definition of $\U$ and $\V$,
note that
$C$ is included to define $\U$
and $A_1$ is excluded to define $\V$,
this is because we want to discuss the problem
in the projective spaces $\PP\mc{S}^n$ and $\PP\U$
whose elements are invariant under scaling.
So we ignored the constraint $A_1\bullet X=1$.
With these flags the optimality
condition~\reff{optcond1}-\reff{optcond2} can be rewritten as
\begin{equation}
\label{betteropt}
X\in \V^\perp \quad \hbox{and} \quad
 Y\in \U \quad \hbox{and} \quad
X \cdot Y   \,=\, 0.
\end{equation}
Our objective is to count the solutions
of this system of polynomial equations.
Notice that if $(X,Y)$ is a solution
and $\lambda, \mu$ are non-zero scalars then
the pair $(\lambda X, \mu Y)$ is also a solution.
What we are counting are equivalence classes
of solutions (\ref{betteropt})
where $(X,Y)$ and $(\lambda X, \mu Y)$ are
regarded as the same point.

Indeed, from now on, for the rest of this paper,
we pass to the usual setting of algebraic geometry:
we complexify each of our linear spaces and we consider
\begin{equation}
\label{projspaces}
\PP \V \subset \PP \U \subset \PP \mc{S}^n \quad \hbox{and} \quad
\PP \U^\perp \subset \PP \V^\perp \subset \PP (\mc{S}^n)^*.
\end{equation}
Each of these six spaces is a {\bf complex projective space}. Note that
\begin{equation}
\label{notethat}
 \dim(\PP \U) = m \quad \hbox{and} \quad
\dim(\PP \V^\perp) =\binom{n+1}2-m.
\end{equation}
If $\W$ is any of the six linear spaces in
\reff{flags} then we write $D^r_\W$ for the
{\em determinantal
variety} of all matrices of rank $\leq r$ in the
corresponding projective space $\PP \W$ in \reff{projspaces}.
Assuming that $D^r_\W$ is non-empty,
the codimension of this variety is independent of
the choice of the ambient space $\PP \W$. By \cite{HT}, we have
 \begin{equation}
 \label{codimformula}
 {\rm codim}(D^r_\W) \quad = \quad
 \binom{n-r+1}{2}.
 \end{equation}
We write $\{XY = 0\}$ for the subvariety of the
product of projective spaces
$\,\PP (\mc{S}^n)^* \times \PP \mc{S}^n\,$
which consists of all pairs $(X,Y)$ of symmetric
$n \times n$-matrices whose matrix product is zero.
If we fix the ranks of the matrices $X$ and $Y$
to be at most $n-r$ and $r$ respectively, then we obtain a subvariety
\begin{equation}
\label{irreduciblesub}
\{XY = 0\}^r  \quad := \quad
 \{XY = 0 \} \,\,\cap \,\, \bigl( D^{n-r}_{(\mc{S}^n)^*} \times D^{r}_{\mc{S}^n} \bigr).
\end{equation}

\begin{lem}
The subvariety $\{XY = 0\}^r$ is irreducible.
\end{lem}

\begin{proof}
Our argument is borrowed from Kempf \cite{GK}.
For a fixed pair of dual bases in $\CC^n$ and $(\CC^n)^*$,
the symmetric matrices $X$ and $Y$ define linear  maps
$$X: \CC^n\to (\CC^n)^*\quad {\rm and} \quad Y: (\CC^n)^*\to \CC^n.$$
Symmetry implies that $\,{\rm ker}(X) \,=\, {\rm Im} (X)^{\perp}\,$
and  ${\rm ker}(Y)\, = \,{\rm Im}(Y)^{\perp}$. Therefore, $\,XY=0\,$
holds if and only if $\,{\rm ker}(X) \supseteq {\rm Im}(Y)$, i.e.
${\rm ker}(Y)\supseteq {\rm ker}(X)^{\perp}$. For fixed rank $r$ and
a fixed subspace $K\subset \CC^n$, the set of pairs $(X,Y)$ such
that $\, {\rm ker} X\supseteq K\,$ and $\, {\rm ker}(Y)\supseteq
K^{\perp}\,$ forms a product of projective linear subspaces of dimension
$\binom{n-r+1}2-1$ and $\binom{r+1}2-1$ respectively.  Over the
Grassmannian $G(r,n)$ of dimension $r$ subspaces in $\CC^n$, these
triples $(X,Y,K)\in \,\PP (\mc{S}^n)^* \times \PP \mc{S}^n \times
G(r,n)$ form a fiber bundle that is an irreducible variety. The
variety  $\{XY = 0\}^r$ is the image under the projection of that
variety onto the first two factors. It is therefore irreducible.
\qed \end{proof}

Thus, purely set-theoretically, our variety has the irreducible decomposition
\begin{equation}
\label{decompo}
 \{XY = 0 \} \quad = \quad \bigcup_{r=1}^{n-1} \{XY = 0\}^r.
\end{equation}
Solving the polynomial equations \reff{betteropt} means
intersecting these subvarieties of
$\,\PP (\mc{S}^n)^* \times \PP \mc{S}^n\,$
with the product of subspaces
 $\,\PP \V^\perp \times \PP \U $.
In other words, what is specified by the
optimality conditions \reff{optcond1}--\reff{optcond2}
is a point in the variety
\begin{equation}
\label{finiteset}
\begin{matrix}
& \,\,\,\, \{XY = 0\}^r \,\,\cap\,\,\,\,\,
\bigl( \PP \V^\perp \times \PP \U \bigr) \,\,\quad \\
\,\, = \,\, &
\{XY = 0 \} \,\,\,\,\cap \,\, \bigl( D^{n-r}_{\V^\perp} \times D^{r}_{\U} \bigr)
\end{matrix}
\end{equation}
which satisfies \reff{optcond3}.
Here $r$ is the rank of the optimal matrix $A(\hat y)$.

In the special case of linear programming,
the matrices $X$ and $Y$ in \reff{finiteset}
are diagonal, and hence
$\{XY = 0\}^r$ is a finite union of linear spaces.
Therefore the optimal pair $(\hat X, \hat Y)$
satisfies some system of linear equations, in accordance
with the fact that the algebraic degree of linear programming
is equal to one.

The following is our main result in this section.

\begin{theorem} \label{zdim}
For generic $\,\U$ and $\V$, the
variety \reff{finiteset} is empty unless
Pataki's inequalities \reff{patiq} hold.
In that case the variety \reff{finiteset}
is reduced, non-empty, zero-dimensional,
and at each point the rank of $X$ and $Y$ is $n-r$ and $r$ respectively.
The cardinality of this variety  depends only on $m,n$ and $r$.
\end{theorem}

The rank condition in Theorem~\ref{zdim} is equivalent to the strict complementarity
condition stated after Theorem~\ref{wsv}.
Theorem~\ref{zdim} therefore implies:

\begin{cor} \label{strictcomp}
The strict complementarity condition holds generically.
\end{cor}

We can now also give an easy proof of Pataki's inequalities:

\smallskip

\noindent 
{\it Proof of Proposition \ref{pataki}}\,
The inequalities \reff{patiq} are implied by the first
sentence of Theorem~\ref{zdim},
since the variety~\reff{finiteset} is not empty.
\qed

\bigskip \noindent
{\it Proof of Theorem \ref{zdim}}\,
Using \reff{notethat} and
 \reff{codimformula}, we
can rewrite \reff{patiq} as follows:
$$
\mbox{codim} (D_{\mc{S}^n}^{r}) \,\leq \, \mbox{dim}
(\PP\U) \qquad \hbox{and} \qquad \mbox{codim} (D_{(\mc{S}^n)^*}^{n-r})
\,\leq\, \mbox{dim} (\PP\V^\perp).
$$
Pataki's inequalities are obviously necessary for
the intersection \reff{finiteset} to be non-empty.
Suppose now that these inequalities are satisfied.

We claim that the dimension of the variety $\{XY = 0\}^r$ equals
$\,\binom{n+1}{2}-2$. In particular, this dimension
is independent of $r$.
To verify this claim, we first note that
there are $\,\binom{n+1}{2} - 1 - \binom{n-r+1}{2}\,$
degrees of freedom in choosing the matrix $\,Y \in D_{\mc{S}^n}^{r}$.
For any fixed matrix $Y$ of rank $r$, consider the
linear system $\,X \cdot Y = 0\,$  of equations for the entries of $X$.
Replacing $Y$ by a diagonal matrix of rank $r$, we see that the
solution space of this linear system has dimension
$\,\binom{n-r+1}{2}-1 $. The sum of these dimensions
equals $\,\binom{n+1}{2}-2\,$ as required. Hence
$$ {\rm dim} \bigl(\{XY = 0\}^r\bigr) \,+\,
{\rm dim} \bigl( \PP \V^\perp \times \PP \U \bigr)
\quad = \quad
\dim( \PP (\mc{S}^n)^{*} \times \PP \mc{S}^n ). $$
By Bertini's theorem \cite[Theorem 17.16]{Har}
for generic choices of the linear spaces $\U$ and $\V$, the
intersection of  $\,\{XY = 0\}^r \,$
with $\,\PP \V^\perp \times \PP \U\,$ is transversal, i.e., in this case finite,
each point of intersection is smooth on both varieties,
and the tangent spaces intersect transversally.
Thus the rank conditions are satisfied at any intersection point.
Furthermore, when the intersection is transversal,
the number of intersection points is independent
of  the choice of $\U$ and $\V$.

We therefore exhibit a non-empty transversal intersection.
Fix any smooth point $(X_{0},Y_{0})\in  D^{n-r}_{({\mc{S}}^n)^{*}} \times D^{r}_{\mc{S}^n}$
with $X_{0} \cdot Y_{0}=0$. This means that $X_{0}$ has rank $n-r$ and $Y_0$ has rank $r$.
After a change of bases, we may assume that
 $X_{0}$ is a diagonal matrix with $r$ zeros and $n \! - \! r$ ones  while $Y_{0}$ has $r$ ones
  and $n \! - \! r$ zeros.

For a generic choice of a subspace
$\V^\perp $ containing  $ X_{0}$
and a subspace  $\U$ containing $ Y_{0}$,
the intersection $\,\{XY = 0\}^r \,\cap \,
(\PP \V^\perp \times \PP \U )\,$ is transversal
away from the point   $\,(X_{0},Y_{0})$.
We must show that the intersection is transversal also at  $(X_{0},Y_{0})$.
We describe the affine tangent space to $\,\{XY = 0\}^r \,$ at  $(X_{0},Y_{0})$
in an affine neighborhood of  $(X_{0},Y_{0})\in \PP (\mc{S}^n)^{*} \times \PP \mc{S}^n$.
The affine neighborhood is the direct sum of the affine spaces parameterized by
\begin{equation}
\label{Xmatrix}
{\bbm
X_{1,1}&\cdots&X_{1,r}&\cdots &\cdots &X_{1,n}\\
\vdots&&&&\ddots& \vdots \\
X_{1,r}&\cdots&X_{r,r}&X_{r,r+1}&\cdots &X_{r,n}\\
X_{1,r+1}&\cdots&X_{r,r+1}&1+X_{r+1,r+1}&\cdots&X_{r+1,n}\\
\vdots&&&&\ddots& \vdots \\
X_{1,n-1}&\cdots&\cdots&\cdots&1+X_{n-1,n-1}&X_{n-1,n}\\
X_{1,n}&\cdots&\cdots&\cdots&X_{n-1,n}&1\\
\ebm} \qquad \hbox{and}
\end{equation}
\begin{equation}
\label{Ymatrix}
{ \bbm
1&Y_{1,2}&\cdots&\cdots &\cdots&Y_{1,n}\\
Y_{1,2}&1+Y_{2,2}&\cdots&\cdots&\cdots &Y_{2,n}\\
\vdots&&&&\ddots & \vdots  \\
Y_{1,r}&\cdots&1+Y_{r,r}&Y_{r,r+1}&\cdots&Y_{r,n}\\
Y_{1,r+1}&\cdots&Y_{r,r+1}&Y_{r+1,r+1}&\cdots&Y_{r+1,n}\\
\vdots&&&&\ddots& \vdots \\
Y_{1,n}&\cdots&\cdots&Y_{r+1,n}&\cdots&Y_{n,n}\\
\ebm.}
\end{equation}
In the coordinates of the matrix equation $XY=0$, the linear terms are
\[
X_{i,j} \hbox{ for } i\leq j\leq r, \,\,\, Y_{i,j} \hbox{ for } r+1\leq i\leq j\, \,\,{\rm and}\,
\,\, X_{i,j}+Y_{i,j} \hbox{ for } i\leq r< j\leq n
\]
These linear forms are independent, and their number equals the
codimension of  $\{XY = 0\}^r$ in $ \PP (\mc{S}^n)^{*} \times \PP \mc{S}^n$.
Hence their vanishing defines the affine tangent space at the smooth point  $(X_{0},Y_{0})$.
For a generic choice these linear terms are independent also of
 $\V$ and $\U^{\perp}$,
which assures the transversality at $(X_{0},Y_{0})$.
\qed

From Theorem~\ref{zdim} we know that the cardinality
of the finite variety~\reff{finiteset} is independent of the choice of generic
$\mathcal{U}$ and $\mathcal{V}$, and it is positive if and only if Pataki's inequalities~\reff{patiq} hold.
We denote this cardinality by $\delta(m,n,r)$.
Our discussion shows that the number  $\delta(m,n,r)$ coincides
with the {\em algebraic degree of SDP}, which was
defined (at the end of Section 2) as the highest degree
of the minimal polynomials of the optimal solution coordinates
$\hat y_i$ and $\hat x_{jk}$.

\section{A table and some formulas}
\setcounter{equation}{0}

The general problem of semidefinite programming can be
formulated in the following primal-dual form which was
derived in Section 3. An instance of SDP
is specified by a flag of linear subspaces
$\,\V \subset \U \subset \mc{S}^n \,$
where ${\rm dim}(\V) = m-1$ and ${\rm dim}(\U) = m+1$,
and the goal is to find
matrices $X,Y \in \mc{S}^n$ such that
\begin{equation}
\label{betteropt2}
X\in \V^\perp \quad \hbox{and} \quad
 Y\in \U \quad \hbox{and} \quad
X \cdot Y   \,=\, 0 \quad \hbox{and} \quad
X,Y \succeq 0.
\end{equation}
Ignoring the inequality constraints
$\,X,Y \succeq 0\,$ and fixing the rank
of $X$ to be $n-r$, the task amounts
to computing the finite variety~\reff{finiteset}.
The algebraic degree of SDP, denoted
$\delta(m,n,r)$, is the cardinality of this projective variety
over $\CC$.
A first result about algebraic degree is the following duality relation.

\begin{prop} \label{degreeduality}
The algebraic degree of SDP satisfies the duality relation
\begin{equation}
\label{dualityrelation}
 \delta \bigl( m,n,r \bigr) \quad = \quad
\delta \bigl(\hbox{$\binom{n+1}{2}$} - m, n, n-r \bigr).
\end{equation}
\end{prop}

\begin{proof}
For generic $C,A_i,b$, by Corollary~\ref{strictcomp}, the
strict complementarity condition  holds.
So in \reff{optcond1}-\reff{optcond3}, if $\mbox{rank}(X^*)=r$,
then $\mbox{rank}(A(y^*) )=n-r$.
But the dual problem~\reff{lmi1}-\reff{lmi2}
can be written equivalently
as some particular primal SDP~\reff{sdp1}-\reff{sdp3}
with $m^\prime = \binom{n+1}{2}-m$ constraints.
The roles of $X$ and $A(y)$ are reversed.
Therefore the duality relation stated in (\ref{dualityrelation}) holds.
\qed
\end{proof}

A census of all values of the algebraic degree
of semidefinite programming for
$m \leq 9$ and $n \leq 6$ is given in
Table~\ref{degtab}. Later, we shall propose
a formula for arbitrary $m$ and $n$. First, let us explain how
Table \ref{degtab} can be constructed.

\begin{table}
\begin{center}
\begin{tabular}{|c|cc|cc|cc|cc|cc|} \hline
& \multicolumn{2}{c|}{$n=2$} &
 \multicolumn{2}{c|}{$n=3$} & \multicolumn{2}{c|}{$n=4$} &
\multicolumn{2}{c|}{$n=5$} & \multicolumn{2}{c|}{$n=6$} \\ \hline
$m$ & $\,\,\,r\,\,\,
 $ & ${\rm degree}$  & $
\,\,\,r\,\,\, $ & ${\rm degree}$  & $ \,\,\,r\,\,\, $ & $
{\rm degree}$
& $\,\,\,r\,\,\, $ & ${\rm degree}$  &
$ \,\,\,r\,\,\, $ & ${\rm degree}$ \\  \hline

1& 1&2&   2&3&   3&4&  4&5&  5&6  \\\hline

2& 1&2&   2&6&   3&12&  4&20&  5&30  \\  \hline

3&  & &  2&4&  3&16&   4&40&  5&80 \\
 &  & &  1&4&  2&10&   3&20&  4&35 \\  \hline

4&  & &  &&  3&8&    4& 40 &  5& 120 \\
 &  & &  1&6&    2&30&   3& 90&  4&210 \\  \hline

5&  & &  &&    &&    4&16&  5&96 \\
 &  & &  1&3&  2&42&   3&207&  4& 672 \\  \hline

6&  & &  &&     &&    &&    5&32 \\
 &  & &   &&   2&30&   3&290&  4&1400 \\
 &  & &   &&   1&8&   2&35&  3&112 \\  \hline

7&  & &  &&  2&10&    3&260&  4&2040 \\
 &  & &   &&   1&16&   2&140&  3&672 \\  \hline

8&  & &  &&  &&    3&140&  4&2100 \\
 &  & &   &&  1&12&   2&260&  3&1992 \\  \hline

9&  & &  &&  &&    3&35&  4&1470 \\
 &  & &   &&  1& 4&   2&290&  3&3812 \\  \hline
\end{tabular}
\end{center}
\smallskip
\caption{The algebraic degree $\delta(m,n,r)$ of semidefinite
programming}
\label{degtab}
\end{table}

Consider the polynomial ring $\Q[X,Y]$ in the
$n(n+1)$ unknowns $x_{ij}$ and $y_{ij}$, and
let $\langle XY \rangle$ be the ideal generated
by the entries of the matrix product $XY$.
The quotient  $R = \Q[X,Y]/\langle XY \rangle$ is
the homogeneous coordinate ring of the variety
$\{XY = 0\}$. For fixed rank $r$ we also
consider the prime ideal $\langle XY \rangle^{\{r\}}$
and the coordinate ring $R^{\{r\}} = \Q[X,Y]/\langle XY \rangle^{\{r\}}$
of the irreducible component  $\,\{XY = 0\}^{r}\,$ in \reff{decompo}.
The rings $R$ and $R^{\{r\}}$ are naturally graded
by the group $\ZZ^2$. The
degrees of the generators are
 ${\rm deg}(x_{ij}) = (1,0)$  and  ${\rm deg}(y_{ij}) = (0,1)$.

A convenient tool for computing and
understanding the columns of Table \ref{degtab}
is the notion of the {\em multidegree} in
combinatorial commutative algebra \cite{MS}.
Recall from \cite[Section 8.5]{MS}
that the multidegree of a $\ZZ^d$-graded affine algebra
is a homogeneous polynomial in $d$ unknowns. Its total degree
is the height of the multihomogeneous presentation ideal.
If $d=2$ then we use the term {\em bidegree} for  the multidegree.
Let $\,\mathcal{C}(R;s,t)\,$
and $\,\mathcal{C}(R^{\{r\}};s,t)\,$ be the
{\em bidegree} of the $\ZZ^2$-graded rings $R$
and $R^{\{r\}}$ respectively. Since
the decomposition \reff{decompo} is  equidimensional,
additivity of the multidegree  \cite[Theorem 8.53]{MS} implies
$$\mathcal{C}(R;s,t)\quad = \quad
\sum_{r=0}^n \mathcal{C}(R^{\{r\}};s,t). $$
The following result establishes the connection to
semidefinite programming:

\begin{theorem} \label{bideg}
The bidegree of the variety $\{XY = 0\}^r$
is the generating function for the
algebraic degree of semidefinite programming.
More precisely,
\begin{equation}
\label{left=right}  \mathcal{C}(R^{\{r\}};s,t) \quad = \quad
 \sum_{m=0}^{\binom{n+1}{2}} \delta(m,n,r) \cdot
 s^{\binom{n+1}{2}-m} \cdot t^m.
\end{equation}
\end{theorem}

\begin{proof}
Fix the ideal $\,I \,=\, \langle XY \rangle^{\{r\}}$,
and let $\,{\rm in}(I)\,$ be its initial monomial ideal
with respect to some term order. The bidegree
$\, \mathcal{C}(R^{\{r\}};s,t) \,$ remains unchanged
if we replace $I$ by ${\rm in}(I)$ in the definition of $R^{\{r\}}$.
The same holds for the right hand side of \reff{left=right}
because we can also define $\delta(m,n,r)$
 by using the initial equations ${\rm in}(I)$
in the place of $\{XY = 0\}^{\{r\}}$
in \reff{finiteset}.  By additivity of
the multidegree, it suffices to consider minimal prime
ideals of ${\rm in}(I)$. These are generated by
subsets of the unknowns $x_{ij}$ and $y_{ij}$.
If such a prime contains $\binom{n+1}{2}-m$ unknowns $x_{ij}$ and
$m$ unknowns $y_{ij}$
then its variety intersects
$\,\PP  \V^\perp \times \PP \U\,$ only if
${\rm dim}(\PP \U) = m$. In this case
the intersection consists of one point.
Hence the contribution to the bidegree
equals $\,s^{\binom{n+1}{2}-m} \cdot t^m \,$ as claimed.
\qed \end{proof}

This formula \reff{left=right} is useful for practical computations because,
in light of the degeneration property \cite[Corollary 8.47]{MS},
the bidegree can be read off from any Gr\"obner basis
for the ideal $\langle XY \rangle^{\{r\}}$. Such a Gr\"obner basis can be computed
for small values of $n$ but no general combinatorial construction
of a Gr\"obner basis is known.
Note that if we set $s=t=1$ in $\mathcal{C}(R^{\{r\}};s,t)$ then we recover
the ordinary $\ZZ$-graded degree of
the ideal $\langle X Y \rangle^{\{r\}}$.

\begin{example} \label{M2}
We examine our varieties $\{X Y = 0\}$
for $n=4$ in {\tt Macaulay~2}:
\begin{verbatim}
R = QQ[x11,x12,x13,x14,x22,x23,x24,x33,x34,x44,
       y11,y12,y13,y14,y22,y23,y24,y33,y34,y44];
X = matrix {{x11, x12, x13, x14},
            {x12, x22, x23, x24},
            {x13, x23, x33, x34},
            {x14, x24, x34, x44}};
Y = matrix {{y11, y12, y13, y14},
            {y12, y22, y23, y24},
            {y13, y23, y33, y34},
            {y14, y24, y34, y44}};
minors(1,X*Y) + minors(2,X) + minors(4,Y); codim oo, degree oo
minors(1,X*Y) + minors(3,X) + minors(3,Y); codim oo, degree oo
minors(1,X*Y) + minors(4,X) + minors(2,Y); codim oo, degree oo
\end{verbatim}
This {\tt Macaulay 2} code verifies that three of the five irreducible components
$\{X Y  = 0\}^{\{r\}}$
in \reff{decompo} all have codimension $10$, and it computes their $\ZZ$-graded degree:
\begin{align*}
 \mathcal{C}(R^{\{3\}};1,1) \quad = & \quad
\sum_{m=1}^4 \delta(m,4,3) \quad = & \quad 4 + 12 + 16 + 8 \quad = & \quad 40 ,\\
\mathcal{C}(R^{\{2\}};1,1) \quad = & \quad
\sum_{m=3}^7 \delta(m,4,2) \quad = & \quad 10 + 30 + 42 + 30 + 10 \quad = & \quad 122.
\end{align*}
The summands are the entries in the $n=4$ column in Table \ref{degtab}.
We compute them by modifying
the {\tt degree} command so that it outputs the bidegree.
\qed
\end{example}

We now come to our main result,
which is a collection of explicit formulas
for many of the entries in Table \ref{degtab},
organized by rows rather than columns.

\begin{theorem}\label{degformula}
The algebraic degree of semidefinite programming,
$\delta(m,n,r)$, is determined by the following formulas
for various special values of $m,n,r$:
\begin{enumerate}
\item If the optimal rank $r$ equals $n-1$ then we have
\[
\delta(m,n,n-1) \,\,\,= \,\,\, 2^{m-1}\binom{n}{m}.
 \]
\item If the optimal rank $r$ equals $n-2$ and  $3\leq m\leq 5$ then we have
\begin{align*}
 \delta(3,n,n-2) \quad & = \quad \binom{n+1}{3} ,\\
\delta(4,n,n-2) \quad & = \quad 6 \binom{n+1}{4}, \\
\delta(5,n,n-2) \quad &= \quad 27\binom{n+1}5+3\binom{n+1}4.
\end{align*}

\item If the optimal rank $r$ equals $n-3$ and  $6\leq m\leq 9$ then we have
\begin{align*}
\delta(6,n,n-3) \,\,\,&= \quad 2 \binom{n+2}{6} \,+\, \binom{n+2}{5}, \\
\delta(7,n,n-3) \,\,\, & = \quad 28\binom {n+3}7-12\binom {n+2}6,       \\
\delta(8,n,n-3) \,\,\, & = \quad  248\binom {n+4}8 - 320\binom {n+3}7 + 84\binom{n+2}6, \\
\delta(9,n,n-3) \,\,\, & = \,\,1794\binom {n \! + \! 5}9-3778\binom {n \! + \! 4}8+ 2436
\binom{n \! + \! 3}7-448 \binom {n \! + \! 2}{6}.
\end{align*}
\end{enumerate}
\end{theorem}

\smallskip

Theorem \ref{degformula} combined with Proposition~\ref{degreeduality} explains all data in
Table \ref{degtab}, except for four special values
which were first found using computer algebra:
\begin{equation}
\label{fourvalues}
\delta(6,6,4) = 1400,\,\,
\delta(7,6,4) = 2040,\,\,
\delta(8,6,4) = 2100,\,\,
\delta(9,6,4) = 1470.
\end{equation}
An independent verification of these numbers will be
presented in Example~\ref{wrapup}. This will illustrate our
general formula which is conjectured to hold arbitrary
values of $m, n$ and $r$. That formula is stated in
Theorem~\ref{smooth} and Conjecture~\ref{conj}.

\smallskip

We note that an explicit and completely general formula for the algebraic degree $\delta(m,n,r)$
was recently found by von Bothmer and Ranestad \cite{vBR}.

\section{ Determinantal varieties and their projective duals}
\setcounter{equation}{0}

In Theorem \ref{bideg}, the algebraic degree of SDP was expressed
in terms of the irreducible components defined by the symmetric matrix equation $X Y = 0$.
In this section we relate this equation to the notion of projective
duality, and we interpret $\delta(m,n,r)$ as the degree of the hypersurface dual
to the variety~$D^r_\U$.

Every (complex) projective space $\PP$ has an associated
dual projective space $\PP^{*}$ whose
points $w$ correspond to hyperplanes $\{w = 0\}$ in $\PP$, and vice versa.
Given any (irreducible) variety $V\subset \PP$, one defines the {\em conormal variety} $CV$
of $V$ to be the (Zariski) closure in $\PP^{*}\times\PP$   of the set of pairs
$(w,v)$ where $\{w=0\}\subset\PP$ is a hyperplane tangent to $V$ at a
smooth point $v\in V$.  The projection of $CV$ in $\PP^{*}$ is the
{\em dual variety} $V^{*}$ to $V$.  The fundamental {\em Biduality Theorem}
states that $CV$ is also the conormal variety of $V^{*}$, and therefore
$V$ is the dual variety to $V^{*}$. For proofs and details see
\cite[\S I.1.3]{GKZ} and \cite[\S 1.3]{Tev}.

In our SDP setting, we write $\PP ({\mc{S}^n})^{*}$ for the projective space dual to $\PP {\mc{S}^n}$.
The conormal variety of the determinantal variety $D^r_{\mc{S}^n}$ is well understood:

\begin{prop} \label{conormal} {\rm \cite[Proposition~I.4.11]{GKZ}}
The irreducible variety
$\,\{XY = 0\}^r \,$ in $\, \PP (\mc{S}^n)^{*}\times \PP {\mc{S}^n}\,$
coincides with the conormal variety of the determinantal variety $D^r_{\mc{S}^n}$ and likewise of
$D^{n-r}_{({\mc{S}^n})^{*}}$.  In particular $D^{n-r}_{({\mc{S}^n})^{*}}$
is the dual variety to $D^r_{\mc{S}^n}$.
\end{prop}

\begin{proof}
Consider a symmetric $n\times n$-matrices $Y$ of rank $r$.
We may assume that $Y$ is diagonal, with the first $r$
entries in the diagonal equal to $1$, and the remaining equal to $0$.
In the affine neighborhood of $Y$, given by $Y_{1,1}\not=0$, the matrices have the form
(\ref{Ymatrix}).
The determinantal subvariety $D^r_{\mc{S}^n}$ intersects this neighborhood in the locus where the
size $r+1$ minors vanish.
The linear parts of these minors specify the matrices $X$ that define
the hyperplanes tangent to $D^r_{\mc{S}^n}$ at $Y$.
The only such minors with a  linear part are those that contain the upper left $r\times r$ submatrix.
Furthermore, their linear part is generated by the entries of the lower right $(n-r)\times (n-r)$ matrix.
But the matrices whose only nonzero entries are in this lower right
submatrix are precisely those that satisfy our matrix equation $\,X \cdot Y =0$.
\qed \end{proof}

\begin{theorem}
\label{MakeConnection} Let $\U$ be a generic $(m+1)$-dimensional
linear subspace of $\mc{S}^n$, and consider the determinantal
variety $D^{r}_\U$ of symmetric matrices of rank at most $r$ in $\PP \U$.
Then its dual variety $(D^{r}_\U)^*$ is a hypersurface if and only if
 Pataki's inequalities \reff{patiq} hold, and, in this case,
 the algebraic degree of semidefinite programming
  coincides with the degree of the dual hypersurface:
\[
\delta(m,n,r) \quad = \quad  \mbox{deg}\,  (D^r_\U)^*.
\]
\end{theorem}

\begin{proof}
Recall the codimension $2$ inclusions $\,\V \subset \U \,$ in $\mc{S}^n$
and  $\,\U^{\perp}\subset\V^{\perp} \,$ in $\,{\mc{S}^n}^*$.
The space of linear forms $\U^{*}$ is naturally identified
with the quotient space ${\mc{S}^n}^{*}/{\U^{\perp}}$,
and hence $\PP ({\mc{S}^n}^{*}/{\U^{\perp}})\, =\,\PP\U^{*}$.
The image of the induced rational map  $\,\PP\V^{\perp}\to
\PP \U^{*} \,$ is the projective line $\,\PP^1 := \,\PP (\V^{\perp}/\U^{\perp}) $.
The points on the line $\PP^1$ correspond to the hyperplanes in $\PP\U$ that contain
the codimension $2$ subspace $\PP\V$.
The map $\,\PP\V^{\perp}\to \PP\U^{*}\,$  induces a map in the first factor
$$\pi \,:\, \PP \V^\perp \times \PP \U
\,\, \to \,\, \PP^1 \times\PP\U
\,\,\, \subset \,\,\,  \PP\U^{*}\times\PP\U.$$
Note that the rational map $\pi$ is only defined outside $\PP \U^\perp \times \PP \U$.

We already know that $\delta(m,n,r)$ is the cardinality of the finite variety
\[
 Z \quad := \quad \{XY = 0\}^r \cap
\bigl( \PP \V^\perp \times \PP \U \bigr) \quad = \quad
\{XY = 0 \} \cap  \bigl( D^{n-r}_{\V^\perp} \times D^{r}_{\U} \bigr).
\]
Since $\U^\perp$ is generic inside $\V^\perp$, none of the points of
 $Z$ lies in $ \PP \U^\perp \times \PP \U $, so we can apply
 $\pi$ to $Z$, and the image $\pi(Z)$ is a finite subset of
 $ \,\PP^1\times\PP\U \, \subset \, \PP\U^{*}\times\PP\U$.

 We next prove that $Z$ and $\pi(Z) $ have the same cardinality.
 By Proposition \ref{conormal}, a point on  $\,\{XY = 0\}^r \,$
  is a pair $(X, Y)$ where the hyperplane $\{X=0\}$ is
 tangent to the determinantal variety $D^{r}_{\mc{S}^n}$
  at the point $Y$. Thus, if $(X_{0}, Y_{0})$ and $(X_{1}, Y_{0})$ are distinct points in $Z$
  that have the same image under $\pi$,
    then both $\{X_{0}=0\}$ and $\{X_{1}=0\}$
  contain the tangent space to $D^r_{\mc{S}^n}$ at $Y_{0}$.
 The  same holds for  $\{sX_{0}+tX_{1}=0\}$ for any $s,t \in \CC$. 
   Hence $Z$ contains  the entire line
   $\,\{ (sX_{0}+tX_{1},Y_{0})\,:\,s,t \in\CC\}$, which is a
   contradiction to $Z$ being finite. Therefore $\pi$
   restricted to $Z$ is one to one, and we conclude
   $\,\# \pi(Z) = \delta(m,n,r)$.

A point $(X,Y)$ in $Z$ represents a hyperplane $\{X=0\}$ in $\PP{\mc{S}^n}$ that contains the subspace
 $\PP\V$ and is tangent to $D^{r}_{\mc{S}^n}$ at a point $Y$ in
 $\PP\U$, i.e. it is tangent to the determinantal variety $D^{r}_{\U}$ at this point.
 Consider the map $\,\{X=0\}\,\mapsto\, \{X=0\}\cap \PP\U\,$ which takes
 hyperplanes in  $\PP{\mc{S}^n}$  to hyperplanes in $\PP\U$. This map
 is precisely the rational map $\,\PP\V^{\perp}\to \PP \U^{*} \,$ defined above.
 The image of $\,(X,Y) \in Z\,$ under $\,\pi\,$ thus represents a
 hyperplane $\,\{X=0\}\cap \PP\U\,$
that contains the codimension $2$ subspace $\PP\V$ and  is tangent to
 $D^{r}_{\U}$ at a point $Y$.
 The hyperplanes in $\PP\U$ that contains $\PP\V$ form  the projective line $\,\PP^1 := \,\PP (\V^{\perp}/\U^{\perp}) $ in $\PP\U^{*}$,
 so $\pi(Z)$ is simply the set of all hyperplanes in that $\PP^1$
which are tangent to $D^{r}_{\U}$.
Equivalently, $\pi(Z)$
 is the intersection of the  dual variety
$(D^{r}_\U)^*$ with a general line in $\PP \U^{*}$. Hence $\pi(Z)$ is nonempty if and only if
$(D^{r}_\U)^*$ is a hypersurface, in which case its cardinality
$ \delta(m,n,r)$ is the degree of that hypersurface.
The first sentence in Theorem \ref{zdim} says that $\, Z \not= \emptyset \,$
if and only if Pataki's inequalities hold.
\qed \end{proof}

The examples in the Introduction serve as an illustration for Theorem \ref{MakeConnection}.
In Example~\ref{ex:vinnikov}, the variety $\,D^{1}_\U \,$ is a cubic Vinnikov curve
and its dual curve $\,(D^{1}_\U)^*$ has degree $\,\delta(2,3,1) = 6$. Geometrically,
a general line meets the (singular) curve in Figure 2
in six points, at least two of which are complex.
A pencil of parallel lines in Figure 1
contains six  that are tangent to the cubic.

In Example \ref{ex:cayley},
the variety $\,D^{2}_\U$ is the Cayley cubic in Figure 3,
and the variety $D^{1}_\U$ consists of its four singular points.
The dual surface $\,(D^{2}_\U)^*\,$ is a Steiner surface of degree $4$,
and the dual surface $(D^{1}_\U)^*$ consists of four planes in
 $\,\PP \U^* \simeq \PP^3 $.
Thinking of $\,D^{2}_\U\,$ as a limit of  smooth cubic surfaces,
we see that each of the
planes in  $(D^{1}_\U)^*$ should be counted with multiplicity two,
as the degree of the surface dual to a smooth cubic surface
is $\, 12 = 4+ 2\cdot 4$.

In general, the degree of the hypersurface $(D^{r}_\U)^*$
depends crucially on the topology and singularities of the
primal variety $D^r_\U$. In what follows, we examine
these features of determinantal varieties, starting with the case $\,\U = \mc{S}^n$.

\begin{prop}
\label{detvar} The determinantal variety $D^r_{\mc{S}^n}$ of symmetric matrices of rank
at most $r$ has codimension $\,\binom{n-r+1}{2}$ and is singular precisely
along the subvariety $D^{r-1}_{\mc{S}^n}$ of matrices of rank at most $r-1$. The degree of
$D^r_{\mc{S}^n}$ equals
\[
 \deg (D^r_{\mc{S}^n}) \quad = \quad
\prod_{j=0}^{n-r-1}{\frac{\binom{n+j}{n-r-j}}{\binom{2j+1}{j}}}.
\]
\end{prop}

\begin{proof}
The codimension,
and the facts that $D_{\mc{S}^n}^r$ is singular along $\,D_{\mc{S}^n}^{r-1} \,$,
appears in \cite[Example~14.16]{Har}.
 The formula for the degree is
\cite[Proposition~12]{HT}.
\qed \end{proof}

When $ D^r_{{\U}}$ is finite, then the above formula suffices to determine our
degree.

\begin{cor} The algebraic degree of semidefinite programming
satisfies
$$  \delta(m,n,r) \,\, = \,\,
\prod_{j=0}^{n-r-1}{\frac {\binom{n+j}{n-r-j}}{\binom{2j+1}{j}}}
 \qquad \hbox{\rm provided} \,\,\,\, m\,=\,
\hbox{$\binom{n-r+1}{2}$}. $$
\end{cor}

\begin{proof}
If $\U\subset{\mc S}^n$ is a generic subspace of dimension $\binom{n-r+1}{2}+1$,
then $\PP\U$ and $ D^r_{{\mc S}^n}$ have complementary dimension in
$\PP {\mc S}^n$, so $\, D^r_{\U} \,=\, \PP\U\cap D^r_{{\mc S}^n}\,$ is finite
and reduced, with cardinality
equal to the degree of $\,D^r_{{\mc S}^n}$.
The dual variety of a finite set in $\PP \U$  is a finite union of
hyperplanes in $\PP \U^*$, one for each point.
\qed \end{proof}

For $m=n=3$ and $r=1$ this formula gives us $\,\delta(3,3,1)=4$.
This is the number of singular points on the Cayley cubic
surface in Example \ref{ex:cayley}.

More generally, whenever $D^r_{\U}$ is smooth
(i.e.~when  $D^{r-1}_{\U}=\emptyset$),
 the degree of the dual hypersurface may be computed by the following {\em Class Formula}.
 For any smooth variety $X\subset \PP^m$
 we write $\chi(X)$ for its {\em Euler number}, which is its
 topological Euler-Poincar\'e characteristic. Likewise, we write
 $\chi(X\cap H)\,$ and $\, \chi(X\cap H\cap H') \,$ for the Euler number
 of the intersection of $X$ with one or two general hyperplanes $H$ and $H'$
 in $\PP^m$ respectively.

\begin{prop}\label{classformula} {\rm (Class Formula, \cite[Theorem 10.1]{Tev})}
If $X$ is any smooth subvariety of $\PP^m$ whose
projective dual variety $X^*$ is a hypersurface, then
$${\rm deg}\,X^*\,\,\,= \,\,\, (-1)^{{\rm dim}(X)} \cdot
\bigl(\,\chi (X) \,-\, 2 \cdot \chi(X\cap H)\,+\,\chi(X\cap H\cap H')\,\bigr)$$
\end{prop}

This formula is best explained in the case when $X$ is a curve,
so that $X\cap H$ is a finite set of points.
The Euler number of $X \cap H$ is its cardinality,
i.e. the degree of $X$.
Furthermore $X\cap H\cap H'=\emptyset$, so the Class Formula reduces to
\begin{equation}
\label{curveclass} {\rm deg} \, X^* \quad = \quad - \, \chi (X)\,\, + \,\,  2 \cdot {\rm deg} \, X.
\end{equation}
To see this, we compute the  Euler number $\chi(X)$ using a general pencil of hyperplanes, namely
those containing the codimension $2$  subspace $H\cap H'\subset\PP^m$.
Precisely $\,\hat d={\rm deg} \, X^*\,$ of the hyperplanes  in this pencil
 are tangent to $X$, and each of
these hyperplanes will be tangent at one point and intersect $X$ in ${\rm deg} \, X-2$ further points.
 So the Euler number is ${\rm deg} \, X-1$ for each of these hyperplane sections.
 The other hyperplane sections all have Euler number   ${\rm deg} \, X$
 and are parameterized by the complement of $\hat d$ points in a $\PP^1$.
 By the multiplicative property of the Euler number, the union of the
 smooth hyperplane sections has Euler number $\,(\chi(\PP^1)-\hat d) \cdot ({\rm deg} \, X)$.
 By  additivity of the Euler number on a disjoint union, we get
the Class Formula (\ref{curveclass}) for curves:
 $\,\chi(X)\,\,=\,\,(\chi(\PP^1)-\hat d \,)\cdot ({\rm deg}\,X)\,+\,
\hat d \cdot ({\rm deg} \, X-1)
\, \,= \,\,2 \cdot {\rm deg} \, X-\hat d$.

The number $\chi(D^r_{\U})$ depends only on $m,n$ and $r$, when $\U$ is general, so we set $\,\chi(m,n,r) \,:=\,\chi(D^r_{\U})$.
When $H$ and $H'$ are general, the varieties $D^r_{\U}\cap H$ and $D^r_{\U}\cap H\cap H'$ are
again general determinantal varieties, consisting of matrices of rank $ \leq r$ in a
codimension $1$ (resp. codimension $2$) subspace of $\U$.
The Class Formula therefore implies the following result.

\begin{cor}\label{class}
Suppose that $\, \binom{n-r+1}2\leq m<\binom{n-r+2}2$. Then we have
    $$\delta(m,n,r) \,\,= \,\, (-1)^{m-\binom{n-r+1}2}\cdot \bigl(\,\chi(m,n,r)\,-2\cdot \chi(m-1,n,r)
    +\chi(m-2,n,r) \bigr). $$
\end{cor}

\begin{proof}
The determinantal variety $D^r_{\U}\subset \PP\U=\PP^m$ for a generic $\U\subset\mc{S}^n$ is
nonempty if and only if $m\geq{\rm codim} D^r_{\U}=\binom{n-r+1}2$, and
it is smooth as soon as $D^{r-1}_{\U}$ is empty, i.e. when
$m<\binom{n-r+2}2$. Therefore the Class Formula applies and gives the
expression for the degree of the dual hypersurface $(D^r_{\U})^*$.
\qed \end{proof}

The duality  in Proposition~\ref{degreeduality} states
$\, \delta \bigl( m,n,r \bigr) \, = \, \delta
\bigl(\hbox{$\binom{n+1}{2}$} - m, n, n-r \bigr)$.
So, as long as one of the two satisfies the inequalities of the
Corollary \ref{class}, the Class Formula computes the  dual
degree.  In Table \ref{degtab}
this covers all cases except $\delta(6,6,4)$, and it covers all
the cases of Theorem \ref{degformula}.  So, for the proof of Theorem
\ref{degformula}, it remains
to compute the Euler number for a smooth $D^r_{\U}$.
We close with the remark that the Class Formula fails
when $D^r_{\U}$ is singular.

\section{Proofs and a conjecture}
\setcounter{equation}{0}

The proof of Theorem \ref{degformula}
has been reduced, by Corollary~\ref{class}, to computing
the Euler number for a smooth degeneracy locus of symmetric matrices.
We begin by explaining the idea of the proof
in the case of the first formula.

\bigskip
\noindent
{\it Proof of Theorem \ref{degformula} (1)}\,
Since $\,\delta \bigl( m,n,n-1 \bigr)
\, = \, \delta \bigl(\binom{n+1}{2} - m, n, 1\bigr)$,
by Proposition \ref{degreeduality}, we may consider
the variety of symmetric rank $1$ matrices, which  is smooth
and coincides with the second Veronese embedding of $\PP^{n-1}$.
The determinantal variety $D^1_{\U}$ is thus
a linear section of this Veronese embedding of $\PP^{n-1}$.
Equivalently, $D^1_{\U}$ is the Veronese image of a
complete intersection of $m-1$ quadrics in
$\PP^{n-1}$.  Our goal is to compute the Euler number of $D^1_{\U}$.

By the Gauss-Bonnet formula, the Euler number of a smooth variety is the degree of the top Chern class of its tangent bundle.
  Let $h$ be the
class of a hyperplane in $\PP^{n-1}$.
The tangent bundle of $D_{\U}^1$ is the quotient of the tangent bundle of $\PP^{n-1}$ restricted to $D_{\U}^1$ and the normal bundle of
$D_{\U}^1$ in $\PP^{n-1}$.
The total Chern class of our determinantal manifold
$D_\U^1$ is therefore the quotient
$$c(D_\U^1)\,\, = \,\, \frac{(1+h)^n}{(1+2h)^{m-1}}, $$
and the top Chern class $c_{n-m}(D_\U^1)$ is the degree $n-m$
term in $c(D_\U^1)$.
 Similarly the top Chern class
 of $D_\U^1\cap H$ and $D_\U^1\cap H\cap H'$ is
the degree $n-m$ terms of $$\frac{(1+h)^n}{(1+2h)^{m}} \cdot 2h \qquad
{\rm and} \qquad
\frac{(1+h)^n}{(1+2h)^{m+1}} \cdot 4h^2,$$ where the last factor indicates
that we evaluate these classes on $D_\U^1$ and use the fact that $H=2h$.
By Proposition \ref{classformula},
the dual degree $\,{\rm deg}(D_\U^1)^*\,$ is obtained by
evaluating $(-1)^{n-m}$ times the degree $n-m$ term in the expression
$$ \frac{(1+h)^n}{(1+2h)^{m-1}} \, -\, 2\frac{(1+h)^n}{(1+2h)^{m}}2h
\,+\,\frac{(1+h)^n}{(1+2h)^{m+1}}4h^2
\,\,\,= \,\,\,\frac{(1+h)^n}{(1+2h)^{m+1}}. $$
That term equals
$\, \binom{n}{m} \cdot \,h^{n-m} $.
Since the degree of our determinantal variety equals
$\,{\rm deg}\,D_\U^1 = \int_{D_\U^1}h^{n-m}=2^{m-1}$,
we conclude
$\,{\rm deg}\,(D_\U^1)^* \, = \, \binom{n}{m} \cdot \,2^{n-m} $.
\qed

In the general case we rely on a formula of Piotr Pragacz \cite{Prag}. He
uses {\em Schur Q-functions} to define an intersection number on $\PP^m$ with
support on the degeneracy locus of a symmetric morphism of vector bundles.
Our symmetric determinantal varieties are special cases of this.
Pragacz' intersection number
does not depend on the smoothness of the degeneracy
locus, but only in the smooth case he shows that
the intersection number equals the Euler number.
By Corollary~\ref{class} we then obtain a formula for $\delta(m,n,r)$ in
the smooth range.

To present Pragacz' formula, we first need to fix our notation for partitions.
 A {\em partition} $\lambda$ is a finite
weakly decreasing
sequence of nonnegative integers $\lambda_{1}\geq\lambda_{2} \geq
\cdots \geq \lambda_{k}\geq 0$. In contrast to the
usual conventions \cite{Mac}, we
distinguish between partitions that differ only by
a string of zeros in the end. For instance, we consider
$(2,1)$ and $ (2,1,0)$ to be distinct partitions.
The {\em length} of a partition $\lambda$ is the number $k$ of its parts,
while the {\em weight} $|\lambda|=\lambda_{1}+\ldots +\lambda_{k}$
is the sum of its parts.
The sum of two partitions is the partition obtained by adding the corresponding $\lambda_{i}$.
A partition is {\em strict} if $\lambda_{i-1} > \lambda_i$
for all $i=1,\ldots,k-1$.
The special partitions $(k, k-1,k-2,\ldots , 1)$ and $(k, k-1,k-2,\ldots ,
1, 0)$
are denoted by $\rho(k)$ and $\rho_{0}(k)$ respectively.

Now, let $E$ be a rank $n$ vector bundle on $\PP^m$ with Chern
roots $x_{1},\ldots ,x_{n}$,
 and let $\phi: E^*\to E$ be a symmetric morphism.
Consider the degeneracy locus $\,D^{r}(\phi)\,$ of points
in $\PP^m$ where the rank of $\phi$ is at most $r$.
For any strict partition $\lambda$ let
$Q_{\lambda}(E)$ be the Schur $Q$-function in the Chern roots
(see \cite{FP,Mac,Prag} for definitions).
Thus $Q_{\lambda}(E)$ is a symmetric polynomial in $x_{1},\ldots ,x_{n}$
of degree equal to the weight of $\lambda$.
Pragacz defines the intersection number
\begin{equation} \label{pin}
e(D^{r}(\phi)\!) \,=\,
\int_{\PP^m}\sum_{\lambda}(-1)^{|\lambda|} \cdot (\!(\lambda+\rho_{0}(n-r-1))\!) \cdot
Q_{(\lambda+\rho(n-r)\!)}(E) \cdot c(\PP^m),
\end{equation}
where the sum is over all partitions $\lambda$ of length
 $n-r$ and weight $|\lambda |\leq m-\binom{n-r+1}{2}$.
Let us carefully describe the ingredients of the formula.
The factor $c(\PP^m)$ is the total Chern class
$\,(1+h)^{m+1}\,$ of $\PP^m$.
For any strict partition $\lambda:=(\lambda_{1}>\lambda_{2} >
\cdots > \lambda_{k}\geq 0)$, the factor $(\!(\lambda)\!)$ is
an integer which is defined as follows.   It depends on $k$ whether or not $\lambda_{k}=0$.
If the number $k$ is $1$ or $2$ then we set $\,(\!(i)\!)
\,:= \,2^i\,$
and
\begin{align*}
(\!(i,j)\!)\,\,\,\,:= \,\,\,\,\binom {i+j}{i}+\binom {i+j}{i-1}+\ldots+\binom
{i+j}{j+1}.
 \end{align*}
In general, when $k$ is even, we set
\begin{align*}
(\!(\lambda_{1},\ldots,\lambda_{k})\!) \,\,\, :=\,\,\,
{\rm  Pfaff}\big( (\!(\lambda_{s},\lambda_{t})\!)\big)_{s<t}.
\end{align*}
Here ``{\rm Pfaff}'' denotes the {\em Pfaffian}
(i.e.~the square root of the determinant) of a skew-symmetric matrix
of even size. Finally, when $k$ is odd, we set
\begin{align*}
(\!(\lambda_{1},\ldots,\lambda_{k})\!) \,\,\,:= \,\,\,
\sum_{j=1}^k (-1)^{j-1}\cdot 2^{\lambda_{j}}\cdot
(\!(\lambda_{1},\ldots,\widehat{\lambda_{j}},\ldots,\lambda_{k})\!).
\end{align*}

\begin{prop}
{\rm (Pragacz, \cite[Prop.~7.13]{Prag})}  \
If the degeneracy locus $\,D^{r}(\phi)\,$ is smooth and of
maximal codimension $\,\binom{n-r+1}2\,$ then
$\,\chi(D^{r}(\phi))=e(D^{r}(\phi))$.
\end{prop}

In our situation, the morphism $\phi$ arises from
the space $\U$ and  is given by a symmetric matrix
whose entries are linear forms on $\PP^m$.
We thus apply the trick, used in \cite{HT} and also in
\cite[Section 6.4]{FP}, of formally writing
$$  E \,\, = \,\, {\Oo_{\PP^{m}}}(\frac h2)
\, \oplus \,\cdots \,\oplus \, {\Oo_{\PP^{m}}}(\frac h2)
\,\, = \,\, n \, { \Oo_{\PP^{m}}}(\frac h2). $$
Here $\phi$ is a map
$\,n{\Oo_{\PP^{m}}}(-\frac h2)\to n{\Oo_{\PP^{m}}}(\frac h2)\,$
and the determinantal variety $D^r_{\U}$ is the locus
of points where this map has rank at most $r$.
The $n$ Chern roots of the split vector bundle $E$ are all equal to
$h/2$. Thus, in  applying Pragacz' formula (\ref{pin}),
we restrict the Schur Q-functions to the diagonal
$\, x_{1}=\ldots=x_{n}= h/2$.

The result of this specialization
of the Schur Q-function is an expression
\begin{equation}
\label{specialSchurQ}
 Q_{\lambda}(E) \,\,\, = \,\,\, b_{\lambda} (n) \cdot h^{|\lambda|},
\end{equation}
where $b_{\lambda}(n)$ is a polynomial in $n$
with $\,b(\ZZ) \subseteq \ZZ[1/2]$.
We multiply (\ref{specialSchurQ}) with the
complementary Chern class of $\PP^m$, which is the expression
$$ c_{|\lambda|}(\PP^m) \,\,\, = \,\,\,
\binom{m+1}{m-\binom{n-r+1}{2}-|\lambda |} \cdot h^{m-|\lambda|}$$
Pragacz' intersection number (\ref{pin}) now evaluates to
\begin{equation}
\label{specialpin}
\sum_{\lambda}(-1)^{|\lambda|} \cdot
(\!(\lambda+\rho_{0}(n-r-1)\!))
\cdot b_{(\lambda+\rho(n-r)\!)}(n) \cdot
\binom{m+1}{m-\binom{n-r+1}{2}-|\lambda
|}.
\end{equation}
We abbreviate the expression  (\ref{specialpin}) by $\,e(m,n,r)$.
Note that $\,e(m,n,r)$ is a polynomial in $n$ for fixed $m$ and $r$.
In the smooth range we can now apply Corollary~\ref{class} to
obtain a formula for
 the degree of the dual variety $(D_{\U}^{r})^*$:
\[
\delta(m,n,r)\,=\,(-1)^{(m-\binom{r+1}{2})} \cdot
\bigl(e(m,n,r)-2 \cdot e(m-1,n,r)+e(m-2,n,r) \bigr).
\]
This yields a formula for the
algebraic degree of semidefinite programming:

\begin{theorem} \label{smooth}
If $\, \binom{n-r+1}2\leq m<\binom{n-r+2}2 \,$ then
$$
\delta(m,n,r)\,\,=\,\,(-1)^d \cdot \sum_{\lambda}(-1)^{|\lambda|}
\cdot (\!(\lambda+\rho_{0}(n-r-1))\!)
\cdot b_{(\lambda+\rho(n-r))}(n)
\cdot \binom{m-1} {d-|\lambda|},
$$
where $\,d \,= \,m-\binom{n-r+1}{2}\,$ is the dimension of
the variety $\,D^r_{\U}\,$ and
the sum is over all partitions  $\lambda$
of length $n-r$ and weight
$\,|\lambda |\leq m-\binom{n-r+1}{2}$.
\end{theorem}

\begin{proof}
It remains only to apply Corollary~\ref{class}.
Comparing the formulas for $e(m,n,r)$, $e(m-1,n,r)$ and
$e(m-2,n,r)$, the difference is in the number of summands and
in the binomial coefficient in the last factor.  But the relation
$$ \binom{m+1}k-2\binom{m}{k-1}+\binom{m-1}{k-2}\,\,= \,\,\binom{m-1}k$$
holds whenever $k>1$,  while
$\binom{m+1}1-2\binom{m}{0}=\binom{m-1}1$ and
$\binom{m+1}0=\binom{m-1}0$. So the formula for
$\delta(m,n,r)$ reduces to the one claimed in the theorem.
\qed \end{proof}

The formula for $\,\delta(m,n,r)\,$ in Theorem \ref{smooth} is
explicit but quite impractical.
To make it more useful, we present a rule
for computing the polynomials $\,b_{(i_{1},\ldots,i_{k})}(n)\,$
for any $i_{1}>\ldots>i_{k}\geq 0$.
First, let $b_{i}(n)$ be the coefficient of $h^i$ in
$$
\frac {(1+ h/2)^n}{(1- h/2)^n} \,\,\, = \quad
b_{0}(n)\,+ \,b_{1}(n) \cdot h \,+ \,\cdots \,+ \, b_{k}(n) \cdot
h^k \,+ \,\cdots $$
The coefficient $b_i = b_i(n)$ is a polynomial of degree $i$ in
the unknown $n$, namely,
\begin{equation}
\label{bformulas1}
b_0 = 1\,, \,\,b_1 = n\,,\,\,
 b_2 = \frac{1}{2}n^2 \,, \,\,
 b_3 = \frac{1}{6}n^3+\frac{1}{12}n \,, \,\,
b_4 = \frac{1}{24}n^4+ \frac{1}{12} n^2 \,,
\end{equation}
\begin{equation}
\label{bformulas2}
 b_5 = \frac{1}{120}n^5+ \frac{1}{24}n^3 + \frac{1}{80}n \,,\,\,
  b_6 = \frac{1}{720}n^6 + \frac{1}{72}n^4+\frac{23}{1440}n^2\,,
  \, \ldots
\end{equation}
We next set $\,b_{(i,0)}(n)=b_{i}(n)\,$ and
$$ b_{(i,j)}(n) \,\, = \,\,
b_{i}(n) \cdot b_{j}(n) \,- \,
2 \cdot \sum_{k=1}^j \, (-1)^{k-1} \cdot
b_{i+k}(n) \cdot b_{j-k}(n). $$
The general formula is now given by distinguishing three cases:
\begin{align*}
b_{(i_{1},\ldots,i_{k})}(n) & ={\rm  Pfaff}(b_{(i_{s},i_{t})})_{s<t}
& \hbox{if $k$ is even}, \\
b_{(i_{1},\ldots,i_{k})}(n) & =b_{(i_{1},\ldots,i_{k},0)}(n)
& \hbox{if $k$ is odd and $i_k > 0$}, \\
b_{(i_{1},\ldots,i_{k-1},0)}(n)& =b_{(i_{1},\ldots,i_{k-1})}(n)
& \hbox{if $k$ is odd and $i_k = 0$}.
\end{align*}

\bigskip
\noindent
{\it Proof of Theorem \ref{degformula} (2),(3)} \,
All seven formulas are gotten by specializing the formula in Theorem \ref{smooth}.
Let us begin with the first one where $m=3$ and $r=n-2$.
 Here $D^{n-2}_{\U}$ is $0$-dimensional
and both the Euler number and the dual degree computes the number
of points in $D^{n-2}_{\U}$.  The formula says
$$
\delta(3,n,n-2)\,\,\,= \,\,\,\sum_{\lambda}(-1)^{|\lambda|}
\cdot (\!(\lambda_1+1,\lambda_{2})\!) \cdot b_{(\lambda_1+2,\lambda_{1}+1)}(n)
\cdot \binom{2}{0-|\lambda |}. $$
The only partition $(\lambda_1,\lambda_2)$ in the
sum is $(0,0)$.  Hence $\, \delta(3,n,n-2)  \,$ equals
$$
(\!(1,0)\!) \cdot b_{(2,1)}(n)\,\,=\,\,
b_{2}(n)b_{1}(n)-2b_{3}(n)\,\,=\,\,
n\frac{n^2}2-2\frac{2n^3+n}{12}
\,\,=\,\,\binom{n+1}3
$$
Next consider the case $m=4$ and $r=n-2$. Here $D^{n-2}_{\U}$ has
dimension $d=1$, and the sum in our formula is over the two
partitions  $\lambda = (0,0)$ and $\lambda = (1,0)$:
\begin{align*}
  \delta(4,n,n-2) & = & &
    - 3 \cdot (\!(1,0)\!) \cdot b_{(2,1)}(n)\, + \,(\!(2,0)\!) \cdot b_{(3,1)}(n) \\
    & = & &
-3 \cdot \bigl(b_{2}(n)b_{1}(n)-2b_{3}(n) \bigr) \, + \,
3 \cdot \bigl(b_{3}(n)b_{1}(n)-2b_{4}(n) \bigr)
\end{align*}
If we substitute (\ref{bformulas1}) into this expression,
then we obtain $\,6\binom{n+1}4\,$ as desired.

The other five cases are similar. We
derive only one more: for $m=8$ and $r = n-3$,
 the sum is over
four partitions $(0,0,0)$, $(1,0,0)$, $(2,0,0)$ and $(1,1,0)$:
\begin{align*}
\delta(8,n,n-3) \,\,\, & = &&
21 \cdot (\!(2,1,0)\!) \cdot b_{321}(n) \,
- \, 7 \cdot (\!(3,1,0)\!) \cdot  b_{421}(n)  \\
&& \, + & \,\,\, 1 \cdot (\!(4,1,0)\!) \cdot b_{521}(n)
\, + \, 1 \cdot (\!(3,2,0)\!) \cdot b_{431}(n).
\end{align*}
After applying the Pfaffian formulas
$$ b_{ijk} \,= \, b_{ij} \cdot b_{k0}
\,-\, b_{ik} b_{j0} \,+ \,b_{jk} b_{i0}
\quad {\rm and} \quad
(\!(i,j,k)\!) = 2^i (\!(j,k)\!) - 2^j (\!(i,k)\!) + 2^k (\!(i,j)\!), $$
we substitute (\ref{bformulas1})--(\ref{bformulas2})
 into this expression and obtain the desired result.
\qed

\begin{exmp}
\label{wrapup}
Consider the four special entries of Table \ref{degtab} which are
listed in (\ref{fourvalues}). Using the duality
(\ref{dualityrelation}) of Proposition~\ref{degreeduality},
we rewrite these values as
$$ \delta(m,6,2)\,\,=\,\,\delta(21-m,6,4) \qquad \hbox{for} \,\,\,\, m =
15,14,13,12. $$
The last three cases satisfy our hypothesis
 $\, 10 = \binom{6-2+1}2\leq 21-m<\binom{6-2+2}2 = 15 $,
 so Theorem~\ref{smooth} applies and furnishes an independent
 proof of the correctness of these values in (\ref{fourvalues}).
The only remaining entry in Table \ref{degtab} is
 $\delta(6,6,4) = 1400$. The formula in
 Theorem \ref{smooth} correctly predicts that value, too.
 \qed
 \end{exmp}

We conjecture that the formula of Theorem~\ref{smooth} holds in the
general singular case.  The formula does indeed make sense without
the smoothness assumption, and it does give the correct number in
all cases that we have checked. Experts in symmetric
function theory might find it an interesting challenge to
verify that the conjectured formula actually satisfies
 the duality relation (\ref{dualityrelation}).

\begin{conj}\label{conj} The formula for
the algebraic degree of semidefinite programming
in Theorem \ref{smooth}
 holds without the restriction in the
range of $m$.
\end{conj}

\medskip

\noindent {\bf Acknowledgements.} We thank Oliver Labs for preparing
 Figures \ref{fig:vinnikov}, \ref{fig:vinnidual} and \ref{fig:cayley}
for us,  and
Hans-Christian Graf von Bothmer for helpful suggestions on varieties
of complexes. We are grateful to the IMA in Minneapolis
for hosting us during our collaboration on this project.
Bernd Sturmfels was partially supported by the
U.S.~National Science Foundation (DMS-0456960).

\end{document}